\documentclass[a4paper]{article}

\usepackage[english]{babel}
\usepackage[utf8x]{inputenc}
\usepackage[T1]{fontenc}

\usepackage{graphicx}
\usepackage{cite}
\usepackage[hang,small]{caption}
\usepackage{subcaption}
\usepackage{color}
\usepackage{amssymb}
\usepackage{bm}
\usepackage{amsthm,amsmath,bm}
\usepackage{color}

\usepackage{tikz}
\usepackage{pgfplots}

\newcommand{\norm}[1]{\left|#1\right|}

\def\clearpage{\hspace{1cm}}

\newcommand{\calT}{\mathcal{T}}

\newcommand{\bbR}{\mathbb{R}}

\newcommand{\bfH}{\bm{h}}
\newcommand{\bfE}{\bm{e}}

\newcommand{\bfJ}{\bm{j}}

\newcommand{\bfn}{\bm{n}}

\newcommand{\bfv}{\bm{v}}

\newcommand{\bfw}{\bm{w}}

\newcommand{\bfx}{\bm{x}}

\DeclareMathOperator{\bfcurl}{{\bf curl}}
\newcommand{\dt}{\partial_t}

\renewcommand{\hat}{\widehat}

\newcommand{\bcD}{\bm{v}_D}
\newcommand{\bcN}{\bm{t}_N}
\newcommand{\bP}{\bm{P}_h}
\newcommand{\bF}{\bm{F}_h}

\newcommand{\bv}{\bm{v}_h}

\newcommand{\bsh}{\widehat{\bm{\sigma}}_h}
\newcommand{\bvh}{\widehat{\bm{v}}_h}

\newcommand{\SndPK}{\mathbf{S}}
\newcommand{\Ohprod}[2]{( #1, #2)_{\calT_{h}}}
\newcommand{\pOhprod}[2]{\langle #1, #2\rangle_{\partial \calT_{h}}}
\newcommand{\pOhGprod}[2]{\langle #1, #2\rangle_{\partial \calT_{h} \backslash \Gamma}}

\newcommand{\GprodD}[2]{\langle #1, #2\rangle_{ \Gamma_D}}
\newcommand{\GprodN}[2]{\langle #1, #2\rangle_{ \Gamma_N}}

\usepackage[a4paper,top=3cm,bottom=2cm,left=3cm,right=3cm,marginparwidth=1.75cm]{geometry}

\usepackage{amsmath}
\usepackage{graphicx}
\usepackage[colorinlistoftodos]{todonotes}
\usepackage[colorlinks=true, allcolors=blue]{hyperref}

\hypersetup{
	colorlinks   = true, 
	urlcolor     = blue, 
	linkcolor    = black, 
	citecolor   = black 
}

\graphicspath{{figures/}}

\title{Hybridized discontinuous Galerkin methods for wave propagation}

\author{P. Fernandez \footnote{MIT Department of Aeronautics and Astronautics, 77 Massachusetts Ave., Cambridge, MA 02139, USA, email: \texttt{pablof@mit.edu}}, \, A. Christophe \footnote{MIT Department of Aeronautics and Astronautics, 77 Massachusetts Ave., Cambridge, MA 02139, USA, email: \texttt{alexchri@mit.edu}}, \, S. Terrana \footnote{MIT Department of Aeronautics and Astronautics, 77 Massachusetts Ave., Cambridge, MA 02139, USA, email: \texttt{terrana@mit.edu}}, \,  N.~C. Nguyen\footnote{MIT Department of Aeronautics and Astronautics, 77 Massachusetts Ave., Cambridge, MA 02139, USA, email: \texttt{cuongng@mit.edu}}, \, and J. Peraire\footnote{MIT Department of Aeronautics and Astronautics, 77 Massachusetts Ave., Cambridge, MA 02139, USA, email: \texttt{peraire@mit.edu}}}

\date{}

\begin{document}

\maketitle

\begin{abstract}

We present the recent development of hybridizable and embedded discontinuous Galerkin (DG) methods for wave propagation problems in fluids, solids, and electromagnetism. In each of these areas, we describe the methods, discuss their main features, display numerical results to illustrate  their performance, and conclude with bibliography notes. The main ingredients in devising these DG methods  are ({\em i}) a local Galerkin projection of the underlying partial differential equations at the element level onto spaces of polynomials of degree $k$  to parametrize the numerical solution in terms of the numerical trace; ({\em ii}) a judicious choice of the numerical flux to provide stability and consistency; and ({\em iii}) a global jump condition that enforces the continuity of the numerical flux to obtain a global system in terms of the numerical trace. These DG methods are termed hybridized DG methods, because they are amenable to hybridization (static condensation) and hence to more efficient implementations. They share many common advantages of DG methods and possess some unique features that make them well-suited to wave propagation problems.

\end{abstract}

\section{Introduction}

Discontinuous Galerkin (DG) methods possess many attractive properties for wave propagation problems. In particular, they are locally conservative, high-order accurate, amenable to complex geometries and unstructured meshes, low dissipative and dispersive, highly parallelizable, and more stable than continuous Galerkin (CG) methods for convection-dominated problems. As a result, DG methods have been widely used in conjunction with explicit time-marching schemes to simulate wave phenomena. Explicit time-integration schemes, however, often become impractical due to the severe time-step size restriction, an issue that is overcome by implicit time-marching schemes. When they are paired with implicit time-marching schemes, DG methods yield a much larger system of equations than CG methods due to the duplication of degrees of freedom along the element faces. The high computational cost and memory footprint make implicit DG methods considerably more expensive than CG methods for a wide variety of applications.


The hybridizable DG (HDG) methods were introduced in \cite{CockburnGopalakrishnanLazarov09} in the framework of steady-state diffusion as part of the effort of devising efficient implicit DG methods for solving elliptic partial differential equations (PDEs). Indeed, the HDG methods guarantee that only the degrees of freedom  of the approximation of the scalar variable on the interelement boundaries are globally coupled, and that the approximate gradient attains optimal order of convergence for elliptic problems \cite{CockburnDongGuzman08,CockburnGopalakrishnanSayas09,CockburnGuzmanWang09SDG}.  The development of the HDG methods was  subsequently extended to a variety of other PDEs: diffusion problems \cite{ChabaudCockburn12,Huynh2013c}, convection-diffusion problems \cite{CockburnDongGuzmanRestelliSacco,nguyen2009d,nguyen2009c,Ueckermann2010},
incompressible flow \cite{Cockburn2011,Cockburn2010e,CockburnSayasHDGStokes,Nguyen2010,Nguyen2011h}, compressible flows \cite{Nguyen2011b,Nguyen2012,Peraire2010,Schutz2013,Fernandez:2018:entropyStable}, continuum mechanics
\cite{Celiker2010,Nguyen2012,Peraire2010,SoonCockburnStolarski09LE}, time-dependent acoustic and elastic wave propagation \cite{CockburnQuennevilleHDGWave,Nguyen2011k}, the Helmholtz equation\cite{Feng2012,Giorgiani2013,GriesmaierMonk}, the time-harmonic Maxwell's equations \cite{Nguyen2011j,Li2014} with the hydrodynamic model \cite{Vidal-Codina2018}, and the time-dependent Maxwell's equations \cite{Christophe2018}. Since the HDG methods inherit many attractive features of DG methods and offer additional advantages in terms of reduced globally coupled degrees of freedom and enhanced accuracy, they have been widely used in conjunction with implicit time-marching schemes to solve time-dependent problems.

Another appealing feature of the HDG methods is that a superconvergent approximation can be computed through a local (and thus inexpensive and highly parallelizable) post-processing step. The superconvergence property cannot be taken for granted since only some combinations of discontinuous finite element spaces and stabilization functions can ensure that property \cite{cockburn2012conditions,cockburn2012superconvergent}. Recently, the theory of $M$-decompositions has provided a simple sufficient condition for the superconvergence. By comparing the dimensions of the space of the approximate trace with the dimensions of the traces of the local volumetric approximations, the $M$-decompositions provide some guidelines to enrich the gradient space such that the superconvergence is ensured. After being presented for diffusion \cite{cockburn2017superconvergenceI,cockburn2017superconvergenceII,cockburn2017superconvergenceIII}, the $M$-decompositions tool has been successfully applied to devise superconvergent HDG methods for Stokes flows \cite{cockburn2017note} and linear elasticity \cite{CockburnFu2018}.

In the setting of wave propagation problems, the HDG methods compare with other finite element methods favorably because they achieve optimal orders of convergence for both the scalar and gradient unknowns and display superconvergence properties \cite{CockburnQuennevilleHDGWave,Feng2012,GriesmaierMonk,Nguyen2011k}. Recently, explicit HDG methods \cite{StanglmeierNguyenPeraireCockburn16} have been introduced for numerically solving the acoustic wave equation. The explicit HDG methods have the same computational cost as other explicit DG methods and provide optimal convergence rates for all the approximate variables. Furthermore, it displays a superconvergence property in agreement with the theoretical results obtained in \cite{CockburnQuennevilleHDGWave}. In spite of the optimal convergence properties, the HDG methods presented in \cite{Nguyen2011k,StanglmeierNguyenPeraireCockburn16} might not be suitable for long-time computations, due to their energy-dissipative characteristics. The dissipative characteristics of HDG for convection-diffusion systems are investigated in \cite{Fernandez:nonModal:2018}. Indeed, it has been observed that dissipative numerical schemes suffer a loss of accuracy for long-time computations, despite their optimal error estimates. Symplectic Hamiltonian HDG methods introduced in \cite{SANCHEZ2017951} are capable of preserving the Hamiltonian structure of the wave equation, while displaying superconvergence properties. Symplectic HDG methods conserve energy and compare favorably with dissipative HDG methods for long-time simulations.

Further extension of the HDG method leads to the introduction of the embedded DG (EDG) method \cite{GuzeyCockburnStolarski07,Peraire2011} and the interior EDG  (IEDG) method \cite{Fernandez2016,Fernandez2017a,Nguyen2015c}. In this paper, we refer to these DG methods as hybridized DG methods, because they are all amenable to hybridization (static condensation) and hence to more efficient implementations. The essential ingredients of hybridized DG methods are ({\em i}) a local Galerkin projection of the underlying PDEs at the element level onto spaces of polynomials of degree $k$ to parametrize the numerical solution in terms of the numerical trace; ({\em ii}) a judicious choice of the numerical flux to provide stability and consistency; and ({\em iii}) a global jump condition that enforces the continuity of the numerical flux to arrive at a global weak formulation in terms of the numerical trace. The only difference among them lies in the definition of the approximation space for the numerical trace. In particular, the numerical trace space of the EDG method is a subset of that of the IEDG method, which in turn is a subset of that of the HDG method.  While the EDG method and the IEDG method do not have superconvergence properties like the HDG method, they yield a smaller system of equations than the HDG method. Indeed, the EDG method has the same degrees of freedom and sparsity pattern as the static condensation of the CG method.  Since the degrees of freedom of the numerical trace on the domain boundary can be eliminated in the IEDG method, IEDG has even less globally coupled unknowns than the EDG method. Thus, the IEDG method is more computationally efficient than both the EDG method and the HDG method.

The remainder of the paper is organized as follows. In Section 2, we introduce preliminary concepts and the notation used throughout the paper. In Section 3, we describe hybridized DG methods for solving the incompressible and compressible Navier-Stokes equations, and present numerical results to demonstrate their performance for a range of flow regimes and wave phenomena. In Section 4, we focus on HDG methods for linear and nonlinear elastodynamics, and show some convergence results for a thin structure. In Section 5, we introduce HDG methods for time-dependent Maxwell's equations with the divergence-free constraint, and present results to verify the convergence and accuracy order. We conclude the paper with our perspectives on future research in Section 6.

\section{\label{s:preliminaries}Preliminaries}

\subsection{Finite element mesh}

Let $T > 0$ be a final time and let $\Omega \subset \mathbb{R}^d$ be an open, connected and bounded physical domain with Lipschitz boundary $\partial \Omega$. We denote by $\mathcal{T}_h$ a collection of disjoint, regular, $p$-th degree curved elements $K$ that partition $\Omega$\footnote{Strictly speaking, the finite element mesh can only partition the problem domain if $\partial \Omega$ is piecewise $p$-th degree polynomial. For simplicity of exposition, and without loss of generality, we assume hereinafter that $\mathcal{T}_h$ actually partitions $\Omega$.}, and set $\partial \mathcal{T}_h := \{ \partial K : K \in \mathcal{T}_h \} $ to be the collection of the boundaries of the elements in $\mathcal{T}_h$. For an element $K$ of the collection $\mathcal{T}_h$, $F= \partial K \cap \partial \Omega$ is a boundary face if its $d-1$ Lebesgue measure is nonzero. For two elements $K^+$ and $K^-$ of $\mathcal{T}_h$, $F=\partial K^{+} \cap \partial K^{-}$ is the interior face between $K^+$ and $K^-$ if its $d-1$ Lebesgue measure is nonzero. We denote by $\mathcal{E}_h^I$ and $\mathcal{E}_h^B$ the set of interior and boundary faces, respectively, and we define $\mathcal{E}_h := \mathcal{E}_h^I \cup \mathcal{E}_h^B$ as the union of interior and boundary faces. Note that, by definition, $\partial \mathcal{T}_h$ and $\mathcal{E}_h$ are different. More precisely, an interior face is counted twice in $\partial \mathcal{T}_h$ but only once in $\mathcal{E}_h$, whereas a boundary face is counted once both in $\partial \mathcal{T}_h$ and $\mathcal{E}_h$.

\subsection{Finite element spaces}

Let $\mathcal{P}_k(D)$ denote the space of polynomials of degree at most $k$ on a domain $D \subset \mathbb{R}^n$, let $L^2(D)$ be the space of Lebesgue square-integrable functions on $D$, and $\mathcal{C}^0(D)$ the space of continuous functions on $D$. Also, let $\bm{\psi}^p_K$ denote the $p$-th degree parametric mapping from the reference element $K_{ref}$ to an element $K \in \mathcal{T}_h$ in the physical domain, and $\bm{\phi}^p_F$ be the $p$-th degree parametric mapping from the reference face $F_{ref}$ to a face $F \in \mathcal{E}_h$ in the physical domain. We then introduce the following discontinuous finite element spaces in $\mathcal{T}_h$,
\begin{equation}
\label{e:spaces1}
\begin{split}
\bm{\mathcal{Q}}_{h}^k  & = \big\{\bm{r} \in [L^2(\mathcal{T}_h)]^{m \times d} \ : \ (\bm{r} \circ \bm{\psi}_K^p ) \in [\mathcal{P}_k(K_{ref})]^{m \times d} \ \ \forall K \in \mathcal{T}_h \big\} ,  \\
\bm{\mathcal{V}}_{h}^k & = \big\{\bm{w} \in [L^2(\mathcal{T}_h)]^m \ : \ (\bm{w} \circ \bm{\psi}_K^p ) \in [\mathcal{P}_k(K_{ref})]^m \ \ \forall K \in \mathcal{T}_h \big\} , \\
{\mathcal{W}}_{h}^k  & = \big\{\psi \in L^2(\mathcal{T}_h) \ : \ (\psi \circ \bm{\psi}_K^p ) \in \mathcal{P}_k(K_{ref}) \ \ \forall K \in \mathcal{T}_h \big\} , \\
\end{split}
\end{equation}
and on the mesh skeleton $\mathcal{E}_{h}$,
\begin{equation}
\label{e:spaces2}
\begin{split}
\bm{\widehat{\mathcal{M}}}_{h}^k  & = \big\{ \bm{\mu} \in [L^2(\mathcal{E}_h)]^m \ 
\ : \ (\bm{\mu} \circ \bm{\phi}^p_F) \in [\mathcal{P}^k(F_{ref})]^m \, \ \forall F \in \mathcal{E}_h \big\} , \\
\bm{\widetilde{\mathcal{M}}}_{h}^k &  = \big\{ \bm{\mu} \in [\mathcal{C}^0(\mathcal{E}_h)]^m \ 
\ : \ (\bm{\mu} \circ \bm{\phi}^p_F) \in [\mathcal{P}^k(F_{ref})]^m \, \ \forall F \in \mathcal{E}_h \big\} , 
\end{split}
\end{equation}
where $m$ is an integer whose particular value depends on the PDE. Note that $\bm{\widehat{\mathcal{M}}}_{h}^k$ consists of functions which are discontinuous at the boundaries of the faces, whereas $\bm{\widetilde{\mathcal{M}}}_{h}^k$ consists of functions that are continuous at the boundaries of the faces. We also denote by $\bm{\mathcal{M}}_{h}^k$ a traced finite element space that satisfies $\bm{\widetilde{\mathcal{M}}}_{h}^k \subseteq \bm{\mathcal{M}}_{h}^k \subseteq \bm{\widehat{\mathcal{M}}}_{h}^k$. In particular, we define
$$ \bm{\mathcal{M}}_{h}^k  = \big\{ \bm{\mu} \in [L^2(\mathcal{E}_h)]^m \ \\
\ : \ (\bm{\mu} \ \circ \ \bm{\phi}^p_F) \in [\mathcal{P}^k(F_{ref})]^m \, \ \forall F \in \mathcal{E}_h , \ \textnormal{and} \ \bm{\mu}|_{\mathcal{E}^{\rm E}_h} \in [\mathcal{C}^0(\mathcal{E}^{\rm E}_h)]^m    \big\} , $$
where $\mathcal{E}^{\rm E}_h$ is a subset of $\mathcal{E}_h$. Note that $\bm{\mathcal{M}}_{h}^k$ consists of functions which are continuous on $\mathcal{E}^{\rm E}_h$ and discontinuous on $\mathcal{E}^{\rm H}_h := \mathcal{E}_h \backslash \mathcal{E}^{\rm E}_h$. Furthermore, if $\mathcal{E}^{\rm E}_h = \emptyset$ then $\bm{\mathcal{M}}_{h}^k = \bm{\widehat{\mathcal{M}}}_{h}^k$, and if $\mathcal{E}^{\rm E}_h = \mathcal{E}_h$ then $\bm{\mathcal{M}}_{h}^k = \bm{\widetilde{\mathcal{M}}}_{h}^k$.

Due to the discontinuous nature of the approximation spaces in \eqref{e:spaces1}, only the degrees of freedom of the approximate trace of the solution on the mesh skeleton $\mathcal{E}_h$, approximated by functions in $\bm{\mathcal{M}}_h^k$, are globally coupled in hybridized DG methods \cite{nguyen2009d,Nguyen2015c}. Hence, different choices of $\mathcal{E}^{\rm E}_h$ lead to different schemes within the hybridized DG family. We briefly discuss three important choices of $\mathcal{E}^{\rm E}_h$. 
The first one is $\mathcal{E}^{\rm E}_h = \emptyset$ and implies $\bm{\mathcal{M}}_{h}^k = \widehat{\bm{\mathcal{M}}}_{h}^k$. This choice corresponds to the hybridizable discontinuous Galerkin (HDG) method \cite{CockburnGopalakrishnanLazarov09}. The second choice is $\mathcal{E}^{\rm E}_h = \mathcal{E}_h$ which implies $\bm{\mathcal{M}}_{h}^k = \widetilde{\bm{\mathcal{M}}}_{h}^k$ and thus enforces the continuity of the approximate trace on {\em all faces}. This choice corresponds to the embedded discontinuous Galerkin (EDG) method introduced in \cite{GuzeyCockburnStolarski07,Peraire2011}. Since $\widetilde{\bm{\mathcal{M}}}_{h}^k \subset \widehat{\bm{\mathcal{M}}}_{h}^k$, the EDG method has fewer globally coupled degrees of freedom that the HDG method. And the third choice of the approximation space $\bm{\mathcal{M}}_{h}^k$ is obtained by setting  $\mathcal{E}^{\rm E}_h = \mathcal{E}^{\rm I}_h$, which implies $\widetilde{\bm{\mathcal{M}}}_{h}^k \subset  \bm{\mathcal{M}}_{h}^k \subset \widehat{\bm{\mathcal{M}}}_{h}^k$, where the inclusions are strict. The resulting approximation space consists of functions which are discontinuous over the union of the boundary faces $\mathcal{E}^{B}_h$ and continuous over the union of the interior faces $\mathcal{E}^{\rm I}_h$. The resulting method has a characteristic of the HDG method on the boundary faces and a characteristic  of the EDG method on interior faces. Because the approximate trace is taken to be continuous only on the interior faces, we shall name this method interior embedded DG (IEDG) method \cite{Fernandez2016,Fernandez2017a,Nguyen2015c} to distinguish it from the EDG method for which the trace is continuous on all faces. We note that the IEDG method enjoys advantages of both the HDG and the EDG methods. First, IEDG inherits the reduced number of global degrees of freedom of EDG. In fact, thanks to the use of face-by-face local polynomial spaces on $\mathcal{E}_h^B$ in the IEDG method, the degrees of freedom of the approximate trace on $\mathcal{E}_h^B$ can be locally eliminated to yield a global matrix system involving only the degrees of freedom of the numerical trace on the interior faces. As a result, the globally coupled unknowns of the IEDG method are even less than those of the EDG method, and IEDG is more efficient than both EDG and HDG. Second, the IEDG scheme enforces the boundary conditions as strongly as the HDG method, thus retaining the boundary condition robustness of HDG. These features make the IEDG method an excellent alternative to the HDG and EDG methods. For additional details on the efficiency and robustness of HDG, EDG and IEDG, the interested reader is referred to \cite{Nguyen2015c}.

It remains to define inner products associated with our finite element spaces. For functions $a$ and $b$ in $L^2(D)$, we denote $(a,b)_D = \int_{D} a b$ if $D$ is a domain in $\mathbb{R}^d$ and $\left\langle a,b\right\rangle_D = \int_{D} a b$ if $D$ is a domain in $\mathbb{R}^{d-1}$. Likewise, for functions $\bm{a}$ and $\bm{b}$ in $[L^2(D)]^m$, we denote $(\bm{a},\bm{b})_D = \int_{D} \bm{a} \cdot \bm{b}$  if $D$ is a domain in $\mathbb{R}^d$ and $\left\langle \bm{a},\bm{b}\right\rangle_D = \int_{D} \bm{a} \cdot \bm{b}$ if $D$ is a domain in $\mathbb{R}^{d-1}$. For functions $\bm{A}$ and $\bm{B}$ in $[L^2(D)]^{m \times d}$, we denote $(\bm{A},\bm{B})_D = \int_{D} \mathrm{tr}(\bm{A}^T \bm{B})$  if $D$ is a domain in $\mathbb{R}^d$ and $\left\langle \bm{A},\bm{B}\right\rangle_D = \int_{D} \mathrm{tr}(\bm{A}^T \bm{B})$ if $D$ is a domain in $\mathbb{R}^{d-1}$, where $\mathrm{tr} \, ( \cdot) $ is the trace operator of a square matrix. We finally introduce the following element inner products
\begin{equation*}
(a,b)_{\mathcal{T}_h} = \sum_{K \in \mathcal{T}_h} (a,b)_K, \qquad (\bm{a},\bm{b})_{\mathcal{T}_h} = \sum_{K \in \mathcal{T}_h} (\bm{a},\bm{b})_K, \qquad (\bm{A},\bm{B})_{\mathcal{T}_h} = \sum_{K \in \mathcal{T}_h} (\bm{A},\bm{B})_K,
\end{equation*}
and face inner products
\begin{equation*}
\left\langle a,b\right\rangle_{\partial \mathcal{T}_h} = \sum_{K \in \mathcal{T}_h} \left\langle a,b\right\rangle_{\partial K}, \qquad \left\langle \bm{a},\bm{b}\right\rangle_{\partial \mathcal{T}_h} = \sum_{K \in \mathcal{T}_h} \left\langle \bm{a},\bm{b}\right\rangle_{\partial K}, \qquad \left\langle \bm{A},\bm{B}\right\rangle_{\partial \mathcal{T}_h} = \sum_{K \in \mathcal{T}_h} \left\langle \bm{A},\bm{B}\right\rangle_{\partial K}.  
\end{equation*}
These notations and definitions are necessary for the remainder of the paper.

\subsection{Time-marching methods}
\label{sec:timemethods}

We describe time-marching methods to integrate in time the following index-1 differential-algebraic equation (DAE) system:
\begin{subequations}
\label{ch5appendx1:eq8}
\begin{alignat}{2}
\bm M \frac{d \bm u}{d t} + \bm f(\bm u, \bm v, t) & = 0 , \qquad t > 0 , \\
\bm g(\bm u, \bm v, t) & = 0 , \qquad t \geq 0 , 
\end{alignat}
\end{subequations}
with initial condition $\bm u(t = 0) = \bm u_0$ and where $\bm M$ is a matrix. The above DAE system will arise from the hybridized DG discretization of time-dependent PDEs in fluids, solids, and electromagnetism. In this context, $\bm M$ is the so-called mass matrix.

\subsubsection{Linear multistep methods}

We denote by $\bm u^n$ an approximation for the function $\bm u(t)$ at discrete time $t^n = n \, \Delta t$, where $\Delta t$ is the time step and $n$ is an integer.  Linear multistep (LM) methods use information from the previous $s$ steps, $\{\bm u^{n+i}\}_{i=0}^{s-1}$, to calculate the solution at the next step $\bm u^{n+s}$.  When we apply a general LM method to the differential part (\ref{ch5appendx1:eq8}a) and treat the algebraic part (\ref{ch5appendx1:eq8}b) implicitly, we arrive at the following algebraic system:
\begin{subequations}
\label{ch5appendx1:eq9}
\begin{alignat}{2}
 \sum_{i=0}^s \left( a_i \bm M \bm u^{n+i} + \Delta t \, b_i \bm f(\bm u^{n+i},\bm v^{n+i}, t^{n+i}) \right)  & =  0, \\
 \label{ch5appendx1:eq9_}
\bm g(\bm u^{n+s}, \bm v^{n+s}, t^{n+s}) & = 0 . 
\end{alignat}
\end{subequations}
The coefficient vectors $\bm a = (a_0,a_1,\ldots,a_s)$ and $\bm b = (b_0,b_1,\ldots,b_s)$ determine the method. If $b_s=0$, the method is called {\em explicit}; otherwise, it is called {\em implicit}. Note we need to solve the system of equations \eqref{ch5appendx1:eq9_} regardless of whether the LM method is explicit or implicit. For this reason and due to their superior stability properties, implicit methods are usually preferred over explicit methods for the temporal integration of DAE systems arising from the spatial hybridized DG discretization of time-dependent PDEs.

Backward difference formula (BDF) schemes are the most popular LM methods for DAE systems. For a BDF scheme with $s$ steps, the system \eqref{ch5appendx1:eq9} becomes
\begin{subequations}
\label{ch5appendx1:eq9b}
\begin{alignat}{2}
a_s\bm M \bm u^{n+s} + \Delta t \, b_s \bm f(\bm u^{n+s},\bm v^{n+s}, t^{n+s}) +  \sum_{i=0}^{s-1} a_i \bm M \bm u^{n+i} & =  &  0, \\
\bm g(\bm u^{n+s}, \bm v^{n+s}, t^{n+s}) & = & 0 .
\end{alignat}
\end{subequations}

\subsubsection{Implicit Runge-Kutta methods}

The coefficients of an $s$-stage Runge-Kutta (RK) method, $a_{ij}, \, b_i, \, c_i , \ 1 \le i,j \le s$, are usually arranged in the form of a {\em Butcher tableau}:
\begin{equation}
\begin{array}{c|cccc}
c_1 & a_{11} & a_{12} & \ldots   & a_{1s} \\
c_2 & a_{21} & a_{22} & \ldots & a_{2s} \\
\vdots & &  \ldots & \vdots & \vdots \\
 c_s & a_{s1} & a_{s2} &  \ldots & a_{ss} \\
\hline
& b_1 & b_2 &  \ldots & b_s 
\end{array}
\end{equation}
For the family of implicit RK (IRK) methods, the RK matrix $a_{ij}$ must be invertible. Let $d_{ij}$ denote the inverse of $a_{ij}$, and let $\bm u^{n,i}$ be the approximation of $\bm u(t)$ at discrete times $t^{n,i} = (t_n + c_i \Delta t) , \ 1 \le i \le s$. The $s$-stage IRK method for the DAE system (\ref{ch5appendx1:eq8}) can be sketched as follows. First, we solve the following $2s$ coupled systems of equations
\begin{subequations}
\label{ch5appendx1:eq10}
\begin{equation}
\label{ch5appendx1:eq10_1}
\begin{split}
\sum_{j=1}^s d_{ij}\bm M \left(\bm u^{n,j} - \bm u^{n} \right) +  \Delta t  \bm f(\bm u^{n,i}, \bm v^{n,i}, t^{n,i})  &=  0, \qquad i = 1,\ldots, s, \\
\bm g(\bm u^{n,i}, \bm v^{n,i}, t^{n,i}) & =  0, \qquad i = 1,\ldots, s ,
\end{split}
\end{equation}
for $(\bm u^{n,i}, \bm v^{n,i})$. Then we compute $\bm u^{n+1}$ from
\begin{alignat}{2}
\label{ch5appendx1:eq10_2}
\bm u^{n+1} & = & \left(1 -  \sum_{i=1}^s b_i \sum_{j=1}^s d_{ij} \right)\bm u^{n} +  \sum_{j=1}^s e_{j} \bm u^{n,j} ,
\end{alignat}
where $e_j =  \sum_{i=1}^s b_i d_{ij}$. Finally, we solve the following system of equations for $\bm v^{n+1}$:
\begin{equation}
\label{ch5appendx1:eq10_3}
\bm g(\bm u^{n+1}, \bm v^{n+1}, t^{n+1})  =  0 .
\end{equation}
\end{subequations}
Note it is possible to advance the system \eqref{ch5appendx1:eq10_1}-\eqref{ch5appendx1:eq10_2} in time without solving \eqref{ch5appendx1:eq10_3}. Hence, we only need to solve \eqref{ch5appendx1:eq10_3} at the particular time steps that we need $\bm v^{n+1}$ for post-processing purposes.

If the RK matrix $a_{ij}$ is a lower-triangular matrix, then the method is called {\em diagonally implicit} RK (DIRK) scheme \cite{Alexander77}. In the case of a DIRK method, each stage of the system \eqref{ch5appendx1:eq10} can be viewed as a BDF step \eqref{ch5appendx1:eq9b} due to the fact the matrix $d_{ij}$ is lower-triangular.

\subsection{Parallel iterative solvers}


We briefly describe the parallel Newton-Krylov-Schwarz method used to solve the (possibly nonlinear) system of algebraic equations that arises from the temporal discretization of the DAE system \eqref{ch5appendx1:eq8} discussed in the previous section. A detailed description of the iterative solver can be found in \cite{Fernandez2016,Fernandez2017a}.


\subsubsection{Nonlinear solver}


To simplify the notation, we shall drop the superscripts that denote the time steps. At any given time step, the nonlinear system of equations (\ref{ch5appendx1:eq9b}) reads as
\begin{subequations}
\label{NLS}
\begin{alignat}{2}
\bm h (\bm u,  \bm v ) & = \bm{0}, \\[1ex]
\bm g(\bm u, \bm v) & = \bm{0},
\end{alignat}
\end{subequations}
where $\bm h$ and $\bm g$ are the discrete nonlinear residuals associated with (\ref{ch5appendx1:eq9b}a) and (\ref{ch5appendx1:eq9b}b), respectively. We solve this nonlinear system using Newton's method. In particular, the linearization of \eqref{NLS} around a given state vector $(\bar{\bm u}, \bar{\bm v})$ yields the following linear system:
\begin{equation}
\left[
\begin{array}{ccc}
\bm A & \bm B \\
\bm C & \bm D
\end{array}
\right]
\left(
\begin{array}{c}
\delta \bm{u} \\
\delta \bm{v}
\end{array}
\right)
=
- \left(
\begin{array}{c}
\bm h (\bar{\bm u}, \bar{\bm v}) \\
\bm g (\bar{\bm u}, \bar{\bm v})
\end{array}
\right) .
\end{equation}
Here the matrices $\bm{A}$, $\bm{B}$, $\bm{C}$, and $\bm{D}$ have entries $A_{ij} = \frac{\partial h_i(\bar{\bm u}, \bar{\bm v})}{\partial u_j}$, $B_{ij} = \frac{\partial h_i(\bar{\bm u}, \bar{\bm v})}{\partial v_j}$, $C_{ij} = \frac{\partial g_i(\bar{\bm u}, \bar{\bm v})}{\partial u_j}$, $D_{ij} = \frac{\partial g_i(\bar{\bm u}, \bar{\bm v})}{\partial v_j}$, respectively.  Since the matrix $\bm{A}$ has block-diagonal structure due to the discontinuous nature of the approximation spaces defined in Section 2.2, $\delta \bm{u}$ can be readily eliminated to obtain a reduced system in terms of $\delta\bm v$ only
\begin{equation}
\label{LS}
\bm{K} \; \delta  \bm v = \bm{r} \ , 
\end{equation}
where $\bm{K}= \bm{D} - \bm{C} \ \bm{A}^{-1} \bm{B}$ and $\bm{r} = \bm{g} -  \bm{C} \  \bm{A}^{-1} \bm{h}$. This is the global system to be solved at every Newton iteration.

To accelerate the convergence of Newton's iterations we compute an initial guess as a solution of a nonlinear least squares problem in which we seek to minimize the norm of the residuals over a subspace. The subspace consists of solutions already computed from the previous time steps. The  Levenberg--Marquardt algorithm is used to solve the nonlinear least squares problem. Further details can be found in \cite{Fernandez2016,Fernandez2017a}.

\subsubsection{Linear solver}

The linear system (\ref{LS}) is solved in parallel using the restarted GMRES method \cite{sasc86} with iterative classical Gram-Schmidt (ICGS) orthogonalization. In order to accelerate convergence, a left preconditioner $\bm{P}^{-1}$ is used and the linear system (\ref{LS}) is replaced by
\begin{equation}
\label{LS1}
\bm{P}^{-1} \bm{K} \; \delta  {\bm{v}} = \bm{P}^{-1} \bm{r} . 
\end{equation}
A restricted additive Schwarz (RAS) \cite{Cai1999} method with $\delta$-level overlap is used as parallel preconditioner. This approach relies on a decomposition of the unknowns in $\delta {\bm{v}}$ among parallel workers; which is performed as described in \cite{Fernandez2017a}. The RAS preconditioner is  defined as
\begin{equation}
\bm{P}^{-1} = \bm{P}_{RAS_{\delta}}^{-1} := \sum_{i=1}^{N} \bm{R}_{i}^{0} \ \bm{K}_{i}^{-1} \ \bm{R}_{i}^{\delta} ,
\end{equation}
where $ \bm{K}_{i} = \bm{R}_{i}^{\delta} \ \bm{K} \ \bm{R}_{i}^{\delta}$ is the subdomain problem, $\bm{R}_{i}^{\beta}$ is the restriction operator onto the subspace associated to the nodes in the $\beta$-level overlap subdomain number $i$, and $N$ denotes the number of subdomains. In our experience, $\delta = 1$ provides the best balance between communication cost and number of GMRES iterations for almost all problems. In practice, we replace $\bm{K}_{i}^{-1}$ by the inverse of the block incomplete LU factorization with zero fill-in, BILU(0), of $\bm{K}_{i}$, that is, $\bm{K}_{i}^{-1} \approx \bm{\widetilde{U}}_{i}^{-1} \ \bm{\widetilde{L}}_{i}^{-1}$. The BILU(0) factorization in each subdomain is performed in conjunction with a Minimum Discarded Fill (MDF) ordering algorithm \cite{Fernandez2017a}.

\section{Wave propagation in fluids}

In this section, we focus on hybridized DG methods for the incompressible and compressible Navier-Stokes equations. Numerical treatment of shock waves using physics-based shock detection and artificial viscosity is described. Numerical results are presented to demonstrate the performance of the methods. The section is ended with bibliography notes.






\subsection{Incompressible Navier-Stokes equations}

\subsubsection{Governing equations}

The unsteady incompressible Navier-Stokes equations for a Newtonian fluid with Dirichlet boundary conditions are given by
\begin{equation}
\label{e:iNS}
\begin{array}{rcll}
\bm{q} - \nabla \bm{v} & = & 0, &  \quad \mbox{in } \Omega \times (0, T) \\
\frac{\partial \bm{v}}{\partial t} -  \nu \ \nabla \cdot \bm{q} +  \nabla p + \nabla \cdot (\bm{v} \otimes \bm{v}) & = & 0 , &  \quad \mbox{in } \Omega \times (0, T) , \\
\nabla \cdot \bm{v} & = & 0, &  \quad \mbox{in } \Omega \times (0, T), \\
\bm{v} & = & \bm{g}, &  \quad \mbox{on } \partial \Omega \times (0, T) , \\
\bm{v} & = & \bm{v}_0, &  \quad \mbox{on } \Omega \times \{t = 0 \} ,
\end{array}
\end{equation}
where $\nu$ denotes the kinematic viscosity of the fluid, $p$ the pressure, $\bm{v} = (v_1, \dots, v_d)$ the velocity vector, $\bm{v}_0$ is the initial velocity field and satisfies the divergence-free condition $\nabla \cdot \bm{v}_0 = 0$ for all $x \in \Omega$, and $\bm{g}$ is the Dirichlet data and satisfies the compatibility condition $\int_{\partial \Omega} \bm{g} \cdot \bm{n} = 0$ for all $t \in (0,T)$. We shall discuss the treatment of other boundary conditions shortly later.

\subsubsection{Formulation}



HDG methods are the only type of hybridized DG method that has been applied to incompressible flows. The HDG method for the unsteady incompressible Navier-Stokes equations \eqref{e:iNS}, as originally proposed in \cite{NguyenPeraireCockburn10AIAAINS,Nguyen2011h}, reads as follows: Find $\big( \bm{q}_h(t) , \bm{v}_h(t) , p_h(t) , \widehat{\bm{v}}_h(t) \big) \in \bm{\mathcal{Q}}_h^k \times \bm{\mathcal{V}}_h^k \times \mathcal{W}_h^k \times \widehat{\bm{\mathcal{M}}}_h^k$ such that
\begin{subequations}
\label{HDGNS}
\begin{alignat}{2}
(\bm{q}_h, \bm{r})_{\mathcal{T}_h} + (\bm{v}_h, \nabla \cdot \bm{r})_{\mathcal{T}_h} -\left\langle \widehat{\bm{v}}_h, \bm{r} \cdot \bm{n} \right\rangle_{\partial \mathcal{T}_h} & =  0, \\
\Big( \frac{d \bm{v}_h}{dt} , \bm{w}\Big)_{\mathcal{T}_h} + \big(\nu \, \bm{q}_h -  p_h \bm{{I}} -  {\bm{v}}_h \otimes {\bm{v}}_h, \nabla  \bm{w}\big)_{\mathcal{T}_h}
+ \left\langle \widehat{\bm{f}}_h, \bm{w}  \right\rangle_{\partial \mathcal{T}_h}  & =  0 , \\
-( \bm{v}_h, \nabla q)_{\mathcal{T}_h} + \left\langle \widehat{\bm{v}}_h  \cdot \bm{n}, q  \right\rangle_{\partial \mathcal{T}_h} & =  0,\\
\label{HDGNS4}
\left\langle \widehat{\bm{f}}_h , \bm{\mu} \right\rangle_{\partial \mathcal{T}_h \backslash \partial \Omega}  + \left\langle \widehat{\bm{v}}_h - \bm{g} , \bm{\mu} \right\rangle_{\partial \Omega} & =  0,  \\
(p_h,1)_{\mathcal{T}_h} & =  0 , \\
\intertext{for all $(\bm{r},\bm{w},q,\bm{\mu}) \in \bm{\mathcal{Q}}_h^k \times \bm{\mathcal{V}}_h^k \times \mathcal{W}_h^k \times \widehat{\bm{\mathcal{M}}}_h^k$ and all $t \in (0,T)$, and}
(\bm{v}_h(t=0) - \bm{v}_0,\bm{w})_{\mathcal{T}_h} & =  0 , \\
\intertext{for all $\bm{w} \in \bm{\mathcal{V}}_h^k$. The integer $m$ in the definition of the spaces $\bm{\mathcal{Q}}_h^k$, $\bm{\mathcal{V}}_h^k$ and $\widehat{\bm{\mathcal{M}}}_h^k$ in Equations \eqref{e:spaces1}$-$\eqref{e:spaces2} is $m=d$ for the incompressible Navier-Stokes equations. Finally, the numerical flux $\widehat{\bm{f}}_h$ is defined as}
\widehat{\bm{f}}_h ( \widehat{\bm{v}}_h , \bm{v}_h , p_h , \bm{q}_h ; \bm{n} ) =  (-\nu \, {\bm{q}}_h  + {p}_h \bm{I}_d +  \widehat{\bm{v}}_h \otimes \widehat{\bm{v}}_h) \cdot \bm{n} + \bm{S} ( \bm{v}_h , \widehat{\bm{v}}_h ; \bm{n} ) \cdot ( \bm{v}_h - \widehat{\bm{v}}_h ) , 
\end{alignat}
\end{subequations}
where $\bm{n}$ is the unit normal vector pointing outwards from the elements, $\bm{I}_d \in \mathbb{R}^{d \times d}$ is the identity matrix, and $\bm{S} \in \mathbb{R}^{d \times d}$ is the so-called stabilization matrix which may depend on $\bm{v}_h$ and $\widehat{\bm{v}}_h$. The stabilization matrix is usually given by $\bm{S} = \bm{S}^{\mathcal{I}} + \bm{S}^{\mathcal{V}}$, where $\bm{S}^{\mathcal{I}}$ is to stabilize the inviscid (convective) operator and $\bm{S}^{\mathcal{V}}$ is to stabilize the viscous (diffusive) operator. These stabilization matrices are typically chosen as $\bm{S}^{\mathcal{I}} = \tau_i \, \bm{I}_d$ and $\bm{S}^{\mathcal{V}} = \tau_v \, \bm{I}_d$, where $\tau_i$ is the inviscid stabilization parameter and $\tau_v$  is the viscous stabilization parameter. Common choices for the former include local $\tau_i = |\widehat{\bm{v}}_h \cdot \bm{n}|$ and global $\tau_i = \sup_{\partial \mathcal{T}_h} |\widehat{\bm{v}}_h\cdot \bm{n}|$ Lax-Friedrichs type approaches, whereas the later is typically defined as $\tau_v = \nu / \ell$ for some characteristic length scale $\ell$ \cite{NguyenPeraireCockburn10AIAAINS,Nguyen2011h} .

\subsubsection{Boundary conditions}

We discuss the numerical treatment of other boundary conditions. In particular, we consider boundary conditions of the form
\begin{equation}
\label{NBC4}
\begin{split}
\bm{v}  & = \bm{g}_D, \quad \mbox{on } \partial \Omega_D \times (0, T), \\
\bm{B}(\bm{q},\bm{v},p) \cdot \bm{n} & = \bm{g}_N, \quad \mbox{on } \partial \Omega_N \times (0, T) ,
\end{split}
\end{equation}
where $\bm{B}$ is a linear boundary operator, and $\partial \Omega_D$ and $\partial \Omega_N$ are such that $\partial \Omega_D \cup \partial \Omega_N = \partial \Omega$ and $\partial \Omega_D \cap \partial \Omega_N = \emptyset$. In order to incorporate these boundary conditions into the HDG discretization, it suffices to replace Equation \eqref{HDGNS4} by
\begin{equation}
\left\langle \widehat{\bm{f}}_h , \bm{\mu} \right\rangle_{\partial \mathcal{T}_h \backslash \partial \Omega}  + \left\langle \widehat{\bm{v}}_h - \bm{g}_D , \bm{\mu} \right\rangle_{\partial \Omega_D} +  \left\langle \widehat{\bm{b}}_h - \bm{g}_N , \bm{\mu} \right\rangle_{\partial \Omega_N} = 0 , 
\end{equation}
where $\widehat{\bm{b}}_h$ is a discretized version of $\bm{B} \cdot \bm{n}$. Some examples of $\bm{B}$ and the corresponding $\widehat{\bm{b}}_h$ are given in Table \ref{t:NBCtable}. Note that $(\bm{q} - \bm{q}^T) \cdot \bm{n} = \bm{\omega} \times \bm{n}$, where $\bm{\omega}$ denotes the vorticity vector, and thus the third and fourth rows in Table~\ref{t:NBCtable} correspond to boundary conditions on the vorticity. Other linear boundary conditions can be treated in a similar manner.

\begin{table}[h]
\begin{center}
\begin{tabular}{|c|c|c|}
\hline
Boundary Condition Type & $\bm{B}$ & $\widehat{\bm{{b}}}_h$ \\
\hline
{Stress} & $ -\nu \, (\bm{q} + \bm{q}^T) + p \bm{{I}}$ &   $(-\nu \, ({\bm{q}}_h + {\bm{q}}_h^T) + {p}_h \bm{{I}}) \cdot \bm{n} + \bm{S} \cdot (\bm{v}_h - \widehat{\bm{v}}_h)$ \\[1ex]
{Viscous stress}$^\ast$ & $ -\nu \, (\bm{q} + \bm{q}^T)$ &   $-\nu \, ({\bm{q}}_h + {\bm{q}}_h^T)  \cdot \bm{n} + \bm{S} \cdot (\bm{v}_h - \widehat{\bm{v}}_h)$ \\[1ex]
Vorticity + pressure & $ -\nu \, (\bm{q} - \bm{q}^T) + p \bm{{I}}$ &   $(-\nu \, ({\bm{q}}_h - {\bm{q}}_h^T) + {p}_h \bm{I}) \cdot \bm{n} + \bm{S} \cdot (\bm{v}_h - \widehat{\bm{v}}_h)$ \\[1ex]
Vorticity$^{\ast,\dagger}$ & $ -\nu \, (\bm{q} - \bm{q}^T)$ &   $-\nu \, ({\bm{q}}_h - {\bm{q}}_h^T) \cdot \bm{n} + \bm{S} \cdot (\bm{v}_h - \widehat{\bm{v}}_h)$ \\[1ex]
Gradient + pressure & $-\nu \, \bm{q} + p \bm{{I}}$ &   $(-\nu \, {\bm{q}}_h + {p}_h  \bm{{I}}) \cdot \bm{n} + \bm{S} \cdot (\bm{v}_h - \widehat{\bm{v}}_h)$ \\[1ex]
Gradient$^\ast$ & $-\nu \, \bm{q}$ &   $-\nu \, {\bm{q}}_h \cdot \bm{n} + \bm{S} \cdot (\bm{v}_h - \widehat{\bm{v}}_h)$ \\
\hline
\end{tabular}
\end{center}
\caption{Examples of other boundary conditions for the incompressible Navier-Stokes equations. The asterisk symbol $^\ast$ indicates that the average pressure condition $(p_h,1)_{\mathcal{T}_h} =  0$ is also imposed. The dagger symbol $\dagger$ indicates that a Dirichlet boundary condition for the normal component of the velocity has also to be provided on $\partial \Omega_N$.}
\label{t:NBCtable}
\end{table}

\subsubsection{Implementation and local post-processing}



The implementation is discussed in \cite{Nguyen2011h}. In short, two different strategies for the Newton-Raphson linearization are proposed in \cite{NguyenPeraireCockburn10AIAAINS,Nguyen2011h}. In the first strategy, the linearized system is hybridized to obtain a reduced linear system involving the degrees of freedom of the approximate velocity and average pressure. The reduced linear system  has a structure of the saddle point problem. In the second strategy, the augmented Lagrangian method developed for the Stokes equations~\cite{Cockburn2010e,Nguyen2010} is used to solve the linearized system. Within each iteration of the augmented Lagrangian method, a linear system involving the degrees of freedom of the approximate velocity only is solved.


The post-processing procedure proposed in \cite{Cockburn2011,Nguyen2011h} can be used to obtain an exactly divergence-free, $\bm{H}({\rm div})$-conforming approximate velocity $\bm{v}_h^\ast$.  This post-processing procedure is local (i.e.\ it is performed at the element level) and thus adds very little to the overall computational cost. Numerical results presented in \cite{Nguyen2011h} show that the approximate pressure, velocity and velocity gradient converge with the optimal order $k+1$ for diffusion-dominated problems with smooth solutions. In such case, the post-processed velocity $\bm{v}_h^\ast$ converges with the order $k+2$ for $k \ge 1$.

\subsection{Compressible Navier-Stokes equations}



\subsubsection{\label{s:cNS_eq}Governing equations}

The unsteady compressible Navier-Stokes equations read as
\begin{equation}
\label{e:NS}
\begin{array}{rcll}
\displaystyle \bm{q} - \nabla \bm{u}  & = & 0,  \quad \mbox{in } \Omega \times (0, T) , \\
\displaystyle \frac{\partial  \bm{u}}{\partial t} +  \nabla  \cdot  \bm{F}(\bm{u},\bm{q})  & = & 0, \quad \mbox{in } \Omega \times (0, T) , \\
\bm{B} (\bm{u} , \bm{q}) & = & 0,  \quad \mbox{on } \partial \Omega \times (0,T) , \\
\bm{u} - \bm{u}_0  & = & 0,  \quad \mbox{on } \Omega \times \{t = 0\} . 
\end{array}
\end{equation}
Here, $\bm{u} = (\rho, \rho v_{j}, \rho E), \ j=1,...,d$ is the $m$-dimensional ($m = d+2$) vector of conserved quantities (i.e.\ density, momentum and total energy), $\bm{u}_0$ is an initial condition, $\bm{B}$ is a boundary operator, and $\bm{F}(\bm{u},\bm{q})$ are the Navier-Stokes fluxes of dimension $m \times d$, given by the inviscid and viscous terms as
\begin{equation}
\label{flux}
\bm{F}(\bm{u},\bm{q}) = \bm{F}^{\mathcal{I}}(\bm{u}) + \bm{F}^{\mathcal{V}}(\bm{u},\bm{q}) = \left( \begin{array}{c}
\rho v_j \\
\rho v_i v_j + \delta_{ij} p \\
 v_j (\rho E + p)
\end{array}
\right) - \left( \begin{array}{c}
0 \\
\tau_{ij}  \\
v_i \tau_{ij} - f_j
\end{array} 
\right) , \qquad i , \, j = 1 , ... , d , 
\end{equation}
where $p$ denotes the thermodynamic pressure, $\tau_{ij}$ the viscous stress tensor, $f_j$ the heat flux, and $\delta_{ij}$ is the Kronecker delta. For a calorically perfect gas in thermodynamic equilibrium, $p = (\gamma - 1) \, \big( \rho E - \rho \, \norm{\bm{v}}^2 / 2 \big)$, where $\gamma = c_p / c_v > 1$ is the ratio of specific heats and in particular $\gamma \approx 1.4$ for air. $c_p$ and $c_v$ are the specific heats at constant pressure and volume, respectively. For a Newtonian fluid with the Fourier's law of heat conduction, the viscous stress tensor and heat flux are given by
\begin{equation}
\label{closures}
\tau_{ij} = \mu \, \bigg( \frac{\partial v_i}{\partial x_j} + \frac{\partial v_j}{\partial x_i} - \frac{2}{3}\frac{\partial v_k}{\partial x_k}\delta_{ij} \bigg) + \beta \, \frac{\partial v_k}{\partial x_k}\delta_{ij} ,  \qquad \qquad \qquad f_j = - \, \kappa \, \frac{\partial T}{\partial x_j} , 
\end{equation}
where $T$ denotes temperature, $\mu$ the dynamic (shear) viscosity, $\beta$ the bulk viscosity, $\kappa = c_p \, \mu / Pr$ the thermal conductivity, and $Pr$ the Prandtl number. In particular, $Pr \approx 0.71$ for air, and additionally $\beta = 0$ under the Stokes' hypothesis.

\subsubsection{Formulation}

The hybridized DG discretization of the unsteady compressible Navier-Stokes equations \eqref{e:NS} reads as follows: Find $\big( \bm{q}_h(t),\bm{u}_h(t), \widehat{\bm{u}}_h(t) \big) \in \bm{\mathcal{Q}}_h^k \times \bm{\mathcal{V}}_h^k \times \bm{\mathcal{M}}_h^k$ such that
\begin{subequations}
\label{IEDG}
\begin{alignat}{2}
\big( \bm{q}_h, \bm{r} \big) _{\mathcal{T}_h} + \big( \bm{u}_h, \nabla \cdot \bm{r} \big)  _{\mathcal{T}_h} -  \big< \widehat{\bm{u}}_h, \bm{r} \cdot \bm{n} \big> _{\partial \mathcal{T}_h}  & =  0, \\
\Big( \frac{\partial \bm{u}_h}{\partial t}, \bm{w} \Big)_{\mathcal{T}_h} - \Big( \bm{F}(\bm{u}_h,\bm q_h) , \nabla \bm{w} \Big) _{\mathcal{T}_h}  +  \left\langle \widehat{\bm{f}}_h(\widehat{\bm{u}}_h,\bm{u}_h,\bm{q}_h), \bm{w} \right\rangle_{\partial \mathcal{T}_h}  & = 0,  \\
\left\langle \widehat{\bm{f}}_h(\widehat{\bm{u}}_h,\bm{u}_h,\bm{q}_h), \bm{\mu} \right\rangle_{\partial \mathcal{T}_h \backslash \partial \Omega} + \left\langle \widehat{\bm{b}}_h(\widehat{\bm{u}}_h,\bm{u}_h,\bm{q}_h), \bm{\mu} \right\rangle_{\partial \Omega} & =  0 , \\
\intertext{for all $(\bm{r},\bm{w},\bm{\mu}) \in \bm{\mathcal{Q}}_h^k \times \bm{\mathcal{V}}_h^k \times \bm{\mathcal{M}}_h^k$ and all $t \in (0,T)$, and}
\big( \bm{u}_h (t = 0) - \bm{u}_0 , \bm{w} \big)_{\mathcal{T}_h} & =  0 , \\
\intertext{for all $\bm{w} \in \bm{\mathcal{V}}^k_h$. The integer $m$ in the definition of the spaces $\bm{\mathcal{Q}}_h^k$, $\bm{\mathcal{V}}_h^k$ and $\bm{\mathcal{M}}_h^k$ in Equations \eqref{e:spaces1}$-$\eqref{e:spaces2} is $m=d+2$ for the compressible Navier-Stokes equations. Finally, the numerical flux $\widehat{\bm{f}}_h$ is defined as}
\widehat{\bm{f}}_h(\widehat{\bm{u}}_h , \bm{u}_h , \bm{q}_h ; \bm{n}) = \bm{F}(\widehat{\bm{u}}_h,\bm{q}_h) \cdot \bm{n} + \bm{S}(\widehat{\bm{u}}_h , \bm{u}_h ; \bm{n}) \cdot (\bm{u}_h - \widehat{\bm{u}}_h) . 
\end{alignat}
\end{subequations}


$\widehat{\bm{b}}_h$ is a boundary flux and its precise definition depends on the type of boundary condition as discussed in Section \ref{s:cNS_BCs}. Like in the incompressible case, the stabilization matrix $\bm{S} \in \mathbb{R}^{m \times m}$ is usually given by the contribution of inviscid and viscous stabilization terms $\bm{S} = \bm{S}^{\mathcal{I}} + \bm{S}^{\mathcal{V}}$. Several choices for the stabilization of the inviscid fluxes have been proposed in \cite{Fernandez2017a,Fernandez2017,Peraire2010}, including local
\begin{equation}
\bm{S}^{\mathcal{I}} = \frac{1}{2} \big( \bm{A}_n (\widehat{\bm{u}}_h) + | \bm{A}_n (\widehat{\bm{u}}_h) | \big) , \qquad \bm{S}^{\mathcal{I}} = | \bm{A}_n (\widehat{\bm{u}}_h) | , \qquad \bm{S}^{\mathcal{I}} = \lambda_{max} (\widehat{\bm{u}}_h) \ \bm{I}_m , 
\end{equation}
and global
\begin{equation}
\bm{S}^{\mathcal{I}} = \Big( \sup_{\partial \mathcal{T}_h} \ \lambda_{max} \big( \widehat{\bm{u}}_h \big) \Big) \, \bm{I}_m
\end{equation}
approaches. Here, $\bm{A}_n = [\partial \bm{F}^{\mathcal{I}} /\partial \bm{u}] \cdot \bm{n}$ is the Jacobian matrix of the inviscid flux normal to the element face, $\lambda_{max}$ denotes the maximum-magnitude eigenvalue of $\bm{A}_n$, $| \, \cdot \, |$ is the generalized absolute value operator, and $\bm{I}_m$ is the $m \times m$ identity matrix. In order to improve stability, smooth surrogates for the operators $( \ \cdot \ + | \cdot |) / 2$ and $| \cdot |$ above are presented in \cite{Fernandez2017}. The following stabilization matrices for the viscous fluxes have been proposed in \cite{Dahm2014a,Fidkowski2016,nguyen2009d,Peraire2011}:
\begin{equation}
\label{Viscous}
\bm{S}^{\mathcal{V}} = \frac{1}{\ell} \frac{\gamma \mu}{Pr} \, \bm{I}_m , \qquad \qquad \bm{S}^{\mathcal{V}} = \frac{1}{\ell}  \, \bm{n}  \cdot \frac{\partial \bm F^{\mathcal{V}}(\widehat{\bm u}_h, \bm q_h)}{\partial \bm q} \cdot \bm{n}, 
\end{equation}
where $\ell$ is either a viscous length scale $\ell_v$ \cite{Dahm2014a,Fidkowski2016}, a global length scale $\ell_g$ \cite{nguyen2009d,Peraire2011} or a characteristic element size $h$. For low Reynolds number flows (in the case $\ell = \{ \ell_v , \ell_g \}$) or for low cell P\'eclet numbers (in the case $\ell = h$), the viscous stabilization plays an important role in the accuracy and stability of the method. Otherwise, it plays a secondary role and is usually dropped.

We note that, for well-resolved simulations, the choice of the stabilization matrix becomes less critical as the polynomial order $k$ increases since the inter-element jumps and numerical dissipation are of order $O(h^{k+1})$ and $O(h^{2(k+1)})$ \cite{Fernandez:2018:entropyStable,williams2018entropy}, respectively, and thus vanish rapidly with increasing $k$. This may not be the case in under-resolved simulations. A comparison of stabilization matrices for under-resolved turbulent flow simulations is presented in \cite{Fernandez2017}. The relationship between $\bm{S}^{\mathcal{I}}$ and the resulting Riemann solver is also discussed in \cite{Fernandez2017}.

\subsubsection{\label{s:cNS_BCs}Boundary conditions}


%
%
%

The definition of the boundary flux $\widehat{\bm{b}}_h$ depends on the type of boundary condition. For example, at the inflow and outflow sections of the domain, we define the boundary flux $\widehat{\bm{b}}_h$ as
\begin{equation}
\widehat{\bm{b}}_h =  \bm{A}^+_n(\widehat{\bm{u}}_h) \cdot (\bm{u}_h - \widehat{\bm{u}}_h) - \bm{A}^-_n(\widehat{\bm{u}}_h) \cdot (\bm{u}_{\partial \Omega} - \widehat{\bm{u}}_h),
\end{equation}
where $\bm{A}^\pm_n = (\bm{A}_n \pm | \bm{A}_n |)/2$ and $\bm{u}_{\partial \Omega}$ is a boundary state. 
%
%
%
%
%
At a solid surface with no slip condition, we extrapolate density and impose zero velocity as follows
\begin{equation}
\widehat{{b}}_{h,1} = {u}_{h,1} - \widehat{{u}}_{h,1} , \qquad \qquad \widehat{{b}}_{h,i} = - \widehat{{u}}_{h,i} \mbox{ for } 2 \le i \le d+1 . 
\end{equation}
The definition of the last component of $\widehat{\bm{b}}_h$ depends on the type of thermal boundary condition. For isothermal walls, for example, we prescribe the temperature $T_w$ as
\begin{equation}
\widehat{{b}}_{hm} = T_w -  \widehat{T}_h(\widehat{\bm{u}}_h) , 
\end{equation}
where the approximate trace of the temperature $\widehat{T}_h(\widehat{\bm{u}}_h)$ is computed from $\widehat{\bm{u}}_h$. For adiabatic walls, we impose zero heat flux as
\begin{equation}
\widehat{{b}}_{h,m} = \widehat{{f}}_{h,m} . 
\end{equation}
Other boundary conditions can be treated in a similar manner.

\subsubsection{Shock capturing method}

For flows involving shocks, we augment the hybridized DG discretization with the physics-based shock capturing method presented in \cite{Fernandez:shockCapturing:2018}. In short, this shock capturing method increases selected fluid viscosities to stabilize and resolve sharp features, such as shock waves and strong thermal and shear gradients, over the smallest distance allowed by the grid resolution. In particular, the bulk viscosity, thermal conductivity and shear viscosity are given by the contribution of the physical $(\beta_f, \kappa_f, \mu_f)$ and artificial $(\beta^*, \kappa^*, \mu^*)$ values, that is,
$$ \beta = \beta_f + \beta^* , \qquad \qquad \kappa = \kappa_f + \kappa^* = \kappa_f + \kappa_1^* + \kappa_2^* , \qquad \qquad \mu = \mu_f + \mu^* . $$
Shock waves, thermal gradients, and shear layers are stabilized by increasing the bulk viscosity, thermal conductivity, and shear viscosity, respectively. Contact discontinuities are stabilized through one or several of these mechanisms, depending on their particular structure. The thermal conductivity is also augmented in hypersonic shock waves through the term $\kappa_1^*$. The artificial viscosities are devised such that the cell P\'eclet number is of order $1$, and in particular are given by
\begin{subequations}
\label{artTransCoef}
\begin{align}
\beta^{*} &= \Phi_{\beta} \bigg[ \rho \ \frac{k_{\beta} \ h_{\beta}}{k} \ \big( \norm{\bm{v}}^2 + c^{*2} \big) ^{1/2} \ \hat{s}_{\beta} \bigg] , \\
\kappa^* &= \kappa_1^* + \kappa_2^* = \Phi_{\beta} \bigg[ \frac{c_p}{Pr_{\beta}^*} \ \bigg( \rho \ \frac{k_{\beta} \ h_{\beta}}{k} \ \big( \norm{\bm{v}}^2 + c^{*2} \big) ^{1/2} \ \hat{s}_{\beta} \bigg) \bigg] + \Phi_{\kappa} \bigg[ \rho \ c_p \ \frac{k_{\kappa} \ h_{\kappa}}{k} \ \big( \norm{\bm{v}}^2 + c^{*2} \big) ^{1/2} \ \hat{s}_{\kappa} \bigg] , \\
\mu^* &= \Phi_{\mu} \bigg[ \rho \ \frac{k_{\mu} \ h_{\mu}}{k} \ \big( \norm{\bm{v}}^2 + c^{*2} \big) ^{1/2} \ \hat{s}_{\mu} \bigg] , 
\end{align}
\end{subequations}
where $c^*$ is the speed of sound at the critical temperature $T^*$, $\Phi_{ \{ \beta, \kappa, \mu \} } \big[ \cdot \big]$ are smoothing operators (not discussed here), $\hat{s}_{ \{ \beta, \kappa, \mu \} }$ are the bulk viscosity, thermal conductivity and shear viscosity sensors (not discussed here), $Pr_{\beta}^*$ is an artificial Prandtl number relating $\beta^*$ and $\kappa_1^*$, $k_{ \{ \beta , \kappa , \mu \} }$ are positive constants of order $1$, and
\begin{subequations}
\label{e:hs}
\begin{align}
h_{\beta} & = h_{ref} \ \frac{|\nabla \rho|}{\big( \nabla \rho ^t \cdot \bm{M}_h^{-1} \cdot \nabla \rho + \epsilon_h \big)^{1/2}} , \\
h_{\kappa} & = h_{ref} \ \frac{|\nabla T|}{\big( \nabla T ^t \cdot \bm{M}_h^{-1} \cdot \nabla T + \epsilon_h \big)^{1/2}} , \\
h_{\mu} & = h_{ref} \ \sigma_{min} (\bm{M}_h) = h_{ref} \ \inf_{|\bm{a}| = 1} \big\{ \bm{a}^t \cdot \bm{M}_h \cdot \bm{a} \big\} , 
\end{align}
\end{subequations}
are the element size in the direction of the density gradient, the temperature gradient and the smallest element size among all possible directions, respectively. In Eq. \eqref{e:hs}, 
$\bm{M}_h$ denotes the metric tensor of the mesh, $h_{ref}$ the reference element size used in the construction of $\bm{M}_h$, and $\epsilon_{h} \sim \epsilon_{m}^2$ is a constant of order machine epsilon squared. The interested reader is referred to \cite{Fernandez:shockCapturing:2018} for additional details on the shock capturing method.

\subsection{Numerical examples}


We present numerical results for several wave phenomena encountered in fluid mechanics, including acoustic waves, shock waves, and the unstable waves responsible for transition to turbulence in a laminar boundary layer. The stabilization matrix is set to $\bm{S} = \lambda_{max} (\widehat{\bm{u}}_h) \ \bm{I}_m$. 
All results are presented in non-dimensional form. $Pr_f = c_p \, \mu_f / \kappa_f = 0.71$, $\beta_f = 0$ and $\gamma = 1.4$ are assumed in all the test problems.

\subsubsection{Inviscid interaction between a strong vortex and a shock wave}

\begin{figure}[b!]
 \centering
 {\includegraphics[width=0.9\textwidth]{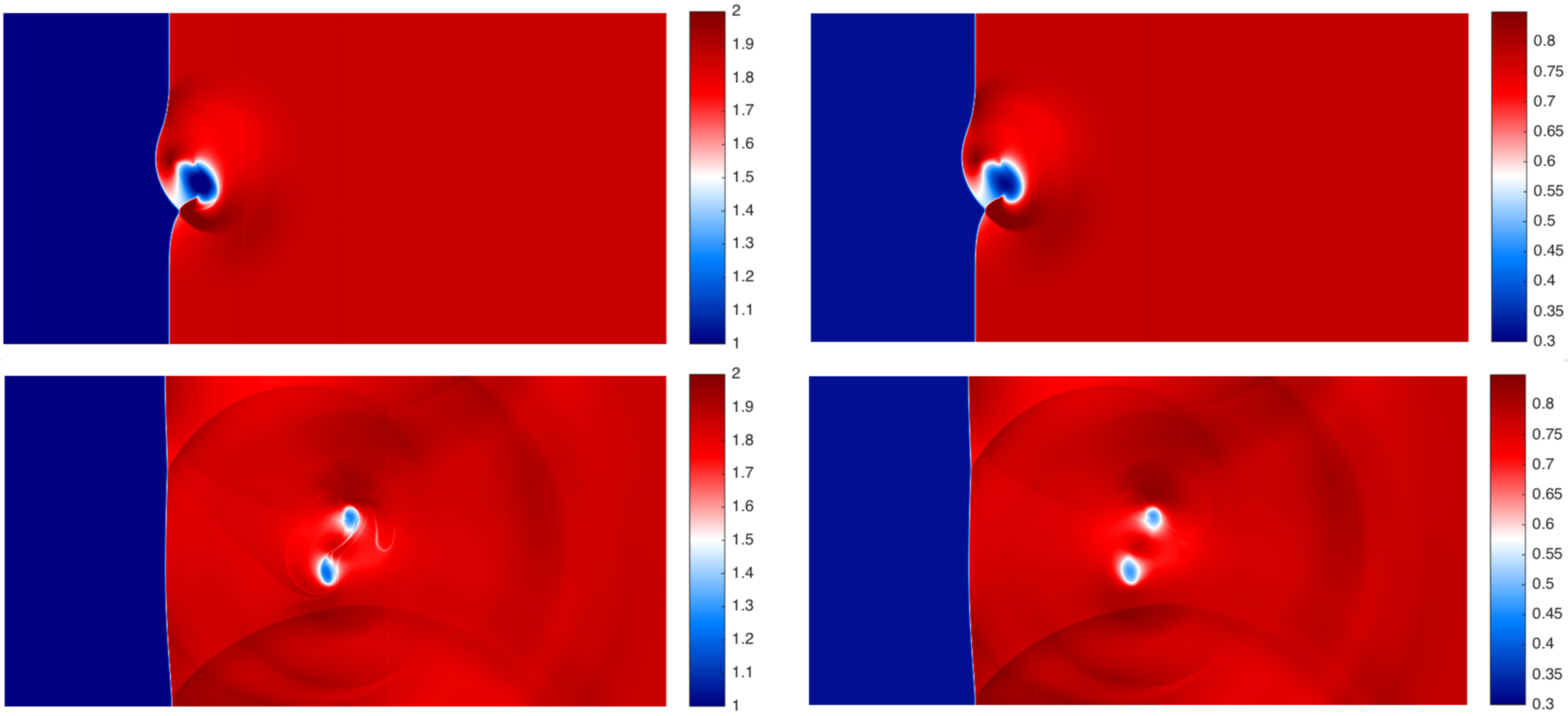}}
  \caption{Non-dimensional density $\rho / \rho_{\infty}$ (left) and pressure $p / (\rho_{\infty} |\bm{v}_{\infty}|^2)$ (right) fields of the strong-vortex/shock-wave interaction problem at the times $t_1$ (top) and $t_2$ (bottom). After the shock wave and the vortex meet, strong acoustic waves are generated and propagate on the downstream side of the shock.}\label{f:densityPressureFields_VSI}
 \end{figure}

\begin{figure}[b!]
 \centering
 {\includegraphics[width=0.9\textwidth]{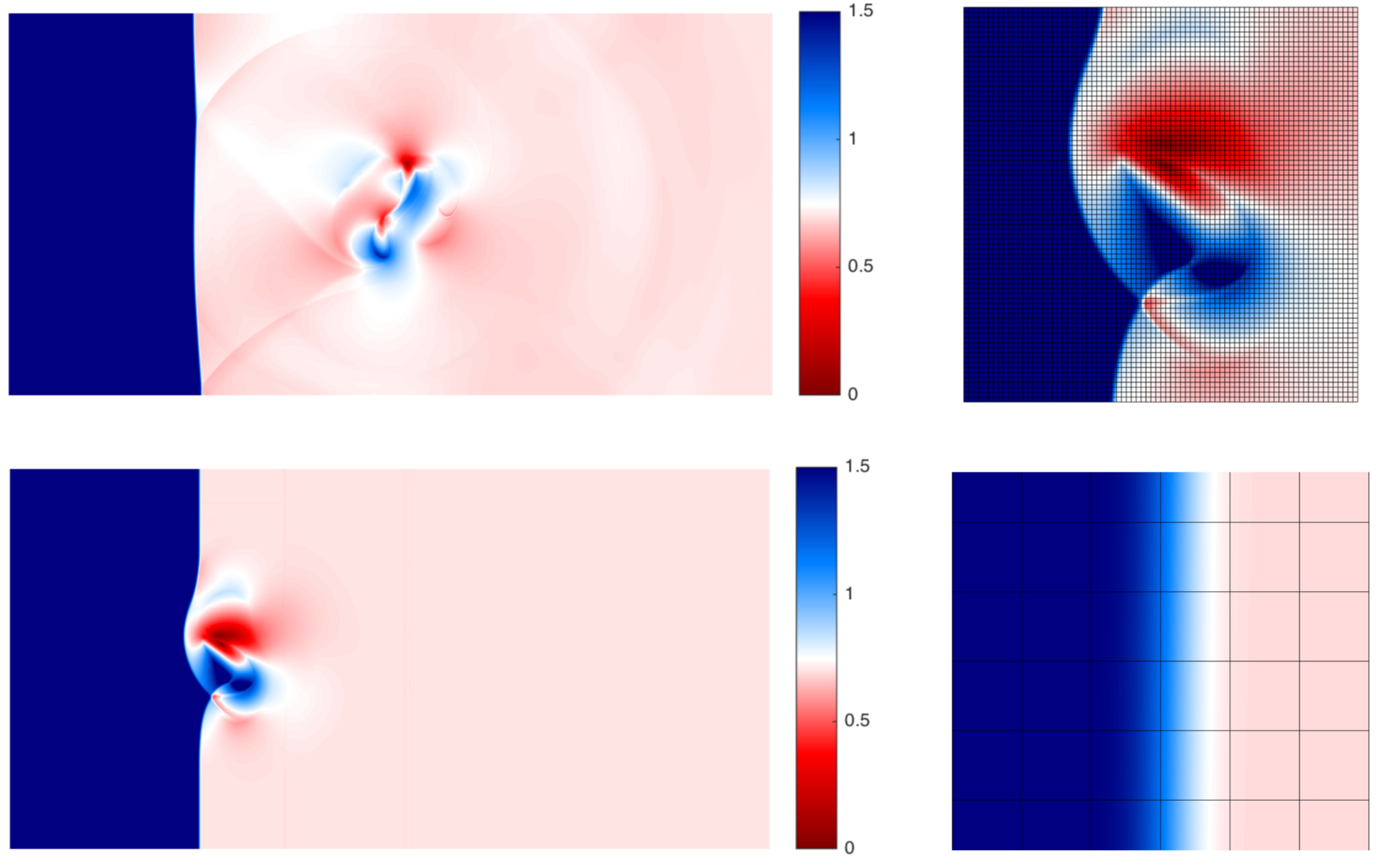}}
 \caption{Mach number field of the strong-vortex/shock-wave interaction problem at the times $t_1$ (top) and $t_2$ (bottom). Zooms around the shock wave are shown on the right images. The shock is non-oscillatory and resolved within one element.}\label{f:machField_VSI}
 \end{figure}

We consider the two-dimensional inviscid interaction between a strong vortex and a shock wave. The problem domain is $\Omega = (0,2 L) \times (0,L)$ and a stationary normal shock wave is located at $x = L / 2$. A counter-clockwise rotating vortex is initially located upstream of the shock and advected downstream by the inflow velocity with Mach number $M_{\infty} = 1.5$. Sixth-order IEDG and third-order DIRK(3,3) schemes are used for the spatial and temporal discretization, respectively. The details of the problem and numerical discretization are presented in \cite{Fernandez2018}. Figure \ref{f:densityPressureFields_VSI} shows the density and pressure fields at the times $t_1 = 0.35 \, \gamma^{1/2} L \, |\bm{v}_{\infty}|^{-1}$ and $t_2 = 1.05 \, \gamma^{1/2} L \, |\bm{v}_{\infty}|^{-1}$. When the shock wave and the vortex meet, the former is distorted and the later split into two separate vortical structures. Strong acoustic waves are then generated from the moving vortex and propagate on the downstream side of the shock. The Mach number fields, together with zooms around the shock wave and the details of the computational mesh, are shown in Figure \ref{f:machField_VSI}. The shock is non-oscillatory and resolved within one element. The shock capturing method does not affect the propagation of the acoustic waves in the sense that it does not introduce artificial dissipation or dispersion \cite{Fernandez2018}.

\subsubsection{Transitional flow over the NACA 65-(18)10 compressor cascade}


We examine the ability of hybridized DG methods to resolve the wave propagation phenomena responsible for natural transition to turbulence in a boundary layer. To this end, we present implicit large-eddy simulation (ILES) results of the three-dimensional NACA 65-(18)10 compressor cascade in design conditions at inlet Reynolds number $Re_1 = 250,000$ and Mach number $M_1 = 0.081$. Third-order IEDG and DIRK(3,3) schemes are used for the discretization. The details of the flow conditions and the numerical setup, as well as the methodology and nomenclature for the boundary layer analysis below, are presented in \cite{Fernandez2017a}.

\begin{figure}[htbp]
\centering
  \includegraphics[width=0.9\textwidth]{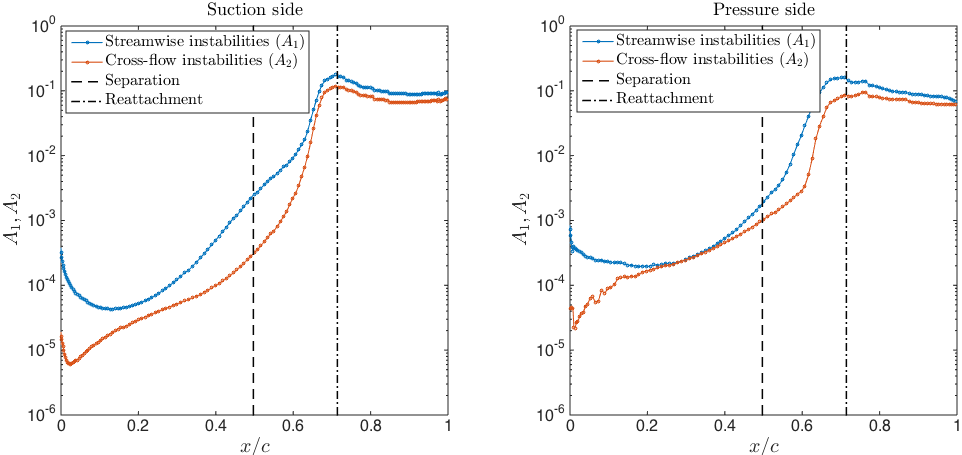}
  \caption{\label{A1A2design} Amplitude of streamwise and cross-flow instabilities on the suction (left) and pressure (right) sides for the NACA 65-(18)10 compressor cascade. The amplitude of the instabilities is non-dimensionalized with respect to the freestream velocity.}
\end{figure}

\begin{figure}
\centering
  \includegraphics[width=0.9\textwidth]{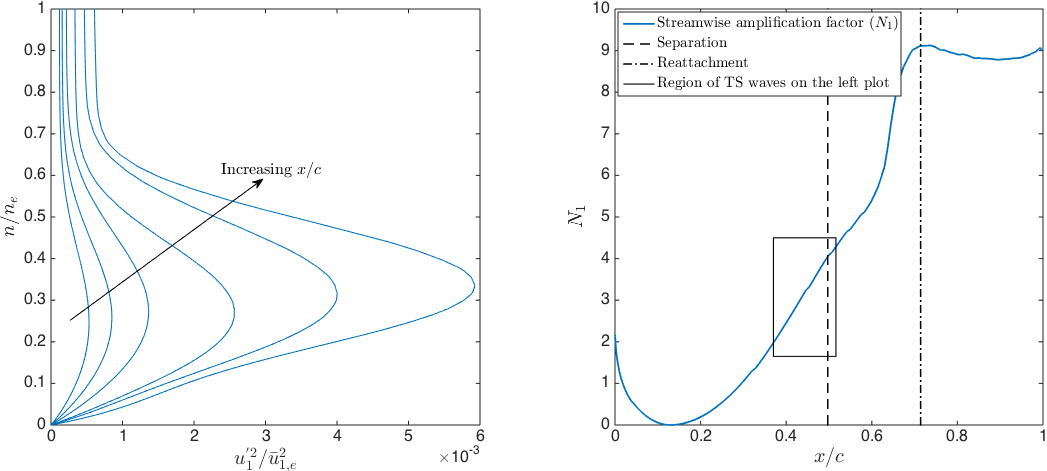}
  \caption{\label{TSwavesSuctionDesign} TS waves (left) and streamwise amplification factor (right) on the suction side for the NACA 65-(18)10 compressor cascade. The box on the right figure indicates the region of the BL in which the TS waves on the left are located.}
\end{figure}

\begin{figure}
\centering
  \includegraphics[width=0.5\textwidth]{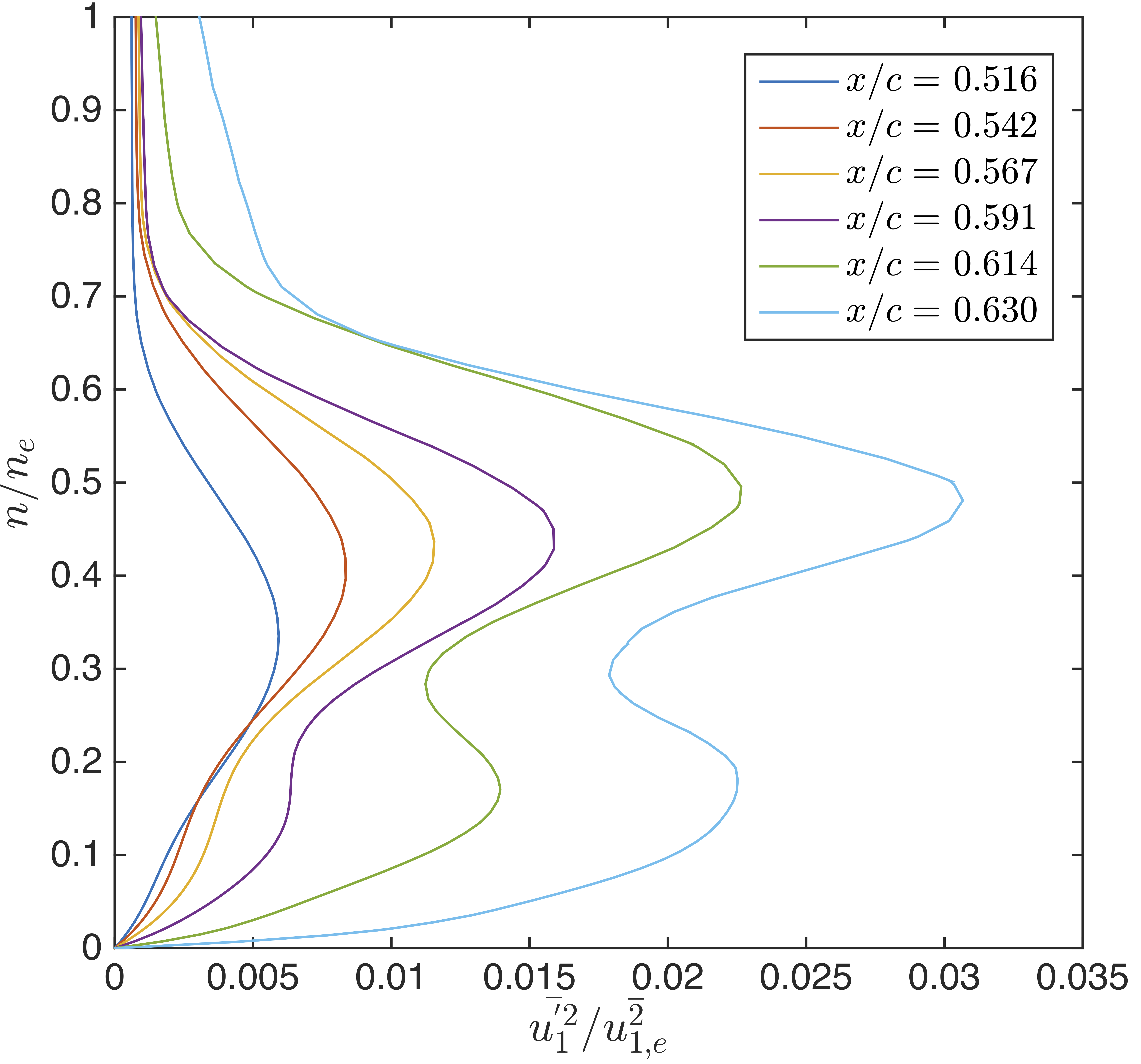}
  \caption{\label{KHinstabilitiesNACAsuction} Transition from TS to KH modes along the separated, suction side boundary layer for the NACA 65-(18)10 compressor cascade.}
\end{figure}

Due to the lack of bypass and forced transition mechanisms and the quasi-2D nature of this flow, natural transition occurs through two-dimensional unstable modes. The two-dimensional nature of transition is illustrated in Figure \ref{A1A2design} through the much larger amplitude of the streamwise instabilities compared to the cross-flow instabilities. In particular, Tollmien-Schlichting (TS) waves form before the boundary layer separates, and Kelvin-Helmholtz (KH) instabilities are ultimately responsible for transition after separation. The former are shown in Figure \ref{TSwavesSuctionDesign} (left) at different BL locations prior to separation. More specifically, the left plot in Figure \ref{TSwavesSuctionDesign} shows the superposition of (1) TS waves and (2) the pressure waves generated in the turbulent boundary layer of the blade at hand and the neighboring blades. The latter effect is responsible for the nonzero fluctuating velocity outside the boundary layer. The growth rate of TS waves along the BL is exponential, as shown on the right of Figure \ref{TSwavesSuctionDesign}\footnote{Note the amplification factor $N_1$ in the $y$-axis is a logarithmic quantity.} and predicted by linear stability theory. It is worth noting the small magnitude of the instabilities compared to the freestream velocity \footnote{Note the amplitude of the instabilities in Figure \ref{A1A2design} is non-dimensionalized with respect to the freestream velocity.}. This shows why very small amount of numerical dissipation is required for transition prediction. Similarly, very low numerical dispersion is needed to properly resolve all the frequencies present in the transition process. After separation, TS waves turn into KH instabilities, as illustrated in Figure \ref{KHinstabilitiesNACAsuction}; which lead to very rapid vortex growth and are ultimately responsible for natural transition in the separated shear layer.


\subsubsection{Transitional flow over the Eppler 387 wing}

We investigate the grid requirements to predict natural transition to turbulence by ILES. In particular, we present grid convergence studies for the transition location of the flow over the three-dimensional Eppler 387 wing at Reynolds numbers of 100,000, 300,000 and 460,000. The Mach number is $M_{\infty} = 0.1$ and the angle of attack $\alpha = 4.0 \ \textnormal{deg}$. Fifth-order HDG and third-order DIRK(3,3) schemes are used for the discretization. Three meshes and non-dimensional time-steps are considered; which correspond to uniform refinement in space and time. The details of these meshes are summarized in Table \ref{eppMeshes}. The interested reader is referred to \cite{Fernandez2016,Fernandez2017a} for additional details on the computational setup.

\begin{table}
\begin{center}
\begin{tabular}{cccccc}
\hline
Mesh No. & $k$ & No. Elements & Element type & Global unknowns & $\Delta t \ |\bm{v}_{\infty}| \, / \, c$ \\
\hline
1 & $4$ & $64,800$ & Tets & $1,959,600 \ \times \ 5$ & $7.937\textnormal{E}-3$ \\
2 & $4$ & $126,360$ & Tets & $3,814,380 \ \times \ 5$ & $6.300\textnormal{E}-3$ \\
3 & $4$ & $254,976$ & Tets & $7,687,680 \ \times \ 5$ & $5.000\textnormal{E}-3$ \\
\hline
\end{tabular}
\end{center}
\caption{\label{eppMeshes}Details of the computational meshes considered for the Eppler 387 wing. $k$ denotes the polynomial order of numerical approximation. $c$ denotes the wing chord. The $\times 5$ factor in {\it global unknowns} accounts for the five components in the Navier-Stokes system.}
\end{table}

\begin{figure}[htbp!]
\centering
\includegraphics[width=1\textwidth]{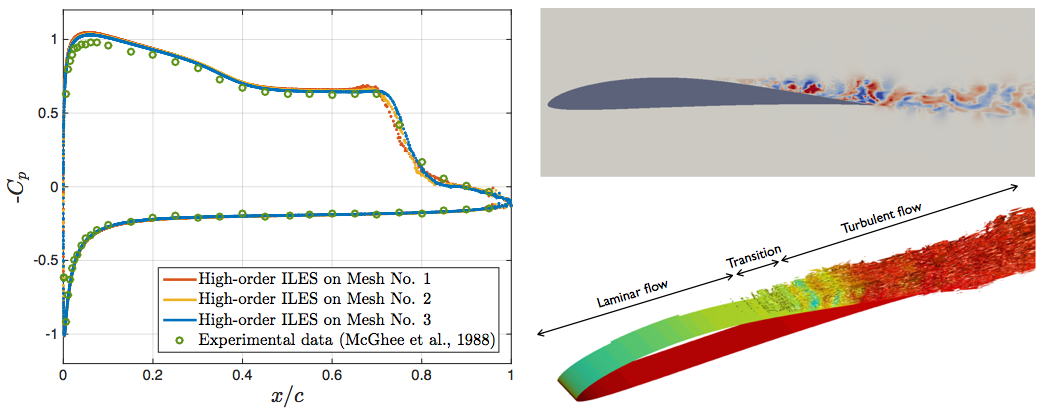}
\caption{\label{CpEppler100k} ILES prediction of the transitional flow over the Eppler 387 wing at $Re = 100,000$: Pressure coefficient (left), instantaneous spanwise velocity (top right), and iso-surface of the Q-criterion colored by pressure (bottom right).}
\end{figure}

\begin{figure}[htbp!]
\centering
\includegraphics[width=1\textwidth]{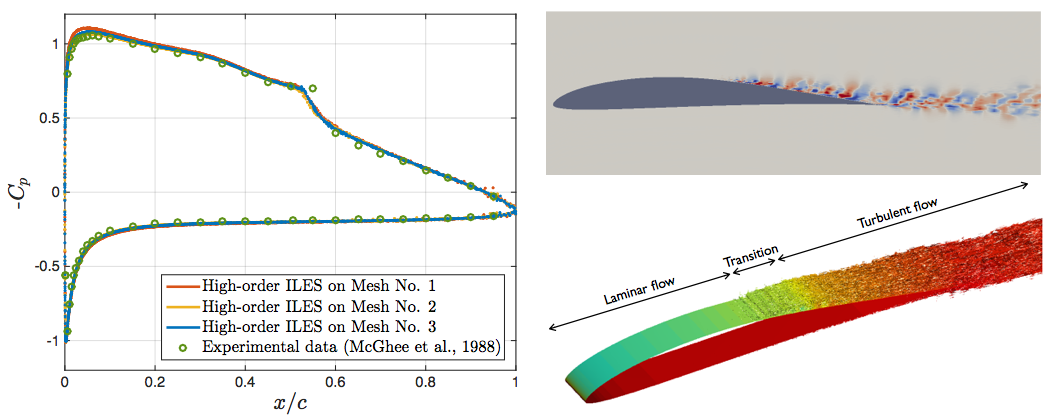}
\caption{\label{CpEppler300k} ILES prediction of the transitional flow over the Eppler 387 wing at $Re = 300,000$: Pressure coefficient (left), instantaneous spanwise velocity (top right), and iso-surface of the Q-criterion colored by pressure (bottom right).}
\end{figure}

\begin{figure}[htbp!]
\centering
\includegraphics[width=1\textwidth]{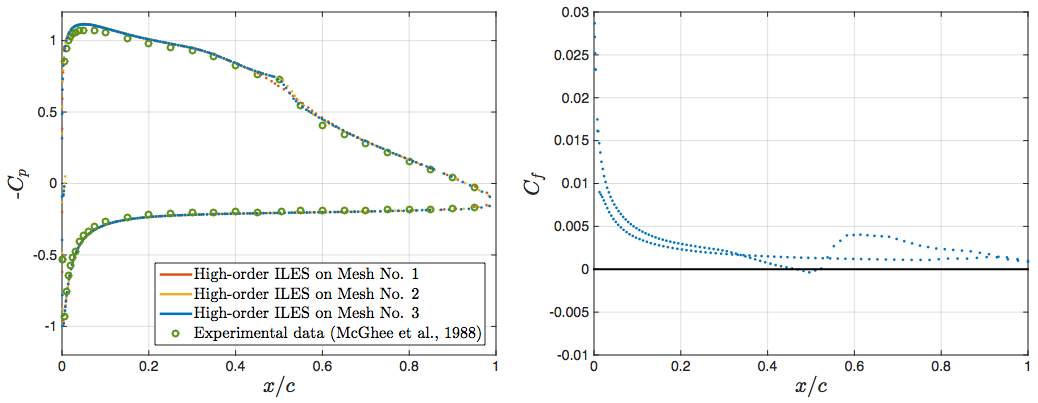}
\caption{\label{CpCfEppler460k} ILES prediction of the pressure coefficient (left) and skin friction coefficient (right) for the Eppler 387 wing at $Re=460,000$.}
\end{figure}

The negative spanwise- and time-averaged pressure coefficient at Reynolds numbers 100,000, 300,000, and 460,000 are shown in Figures \ref{CpEppler100k}, \ref{CpEppler300k} and \ref{CpCfEppler460k}, respectively. The simulation results converge to the experimental data \cite{McGhee1988_2} as the mesh is refined\footnote{The mismatch between the simulation and the experimental data near the leading edge is due to the missing vortex upwash induced by the finite extent of the computational domain, and not due to discretization errors \cite{Fernandez2017a,Fernandez2017}.}. In particular, the error in the transition location is below $0.01c$, $0.005c$, and $0.01c$ at Reynolds number 100,000, 300,000, and 460,000, respectively, even with mesh No. 1. The effective resolution of this mesh is equivalent to a cell-centered finite volume discretization with 691,200 elements. These grid requirements are much smaller than those typically needed with low-order schemes.

The numerical results for the NACA 65-(18)10 cascade and the Eppler 387 wing demonstrate the advantage of high-order DG methods to simulate transitional flows, as they require much fewer elements and degrees of freedom to accurately predict transition than low-order methods. This is justified by the following observation \cite{Fernandez2017a}:  {\it Simulating transition is challenging mostly due to the small magnitude of the instabilities involved, rather than due to their length and time scales}. A low-order scheme may kill the small instabilities because of high numerical dissipation even when the mesh size and time-step size are sufficiently small to represent the length and time scales of the instabilities. We note, however, that high-order methods become more and more computationally expensive (per degree of freedom) as the order of accuracy increases. As discussed in \cite{Fernandez2017a,Fernandez2017}, the hybridized DG methods seem to yield the best trade-off between accuracy and computational cost for transitional flows when the accuracy order is between 3 and 5.

\subsubsection{Transonic flow over the T106C low-pressure turbine}

\begin{figure}[b!]
\centering
\includegraphics[width=1.0\textwidth]{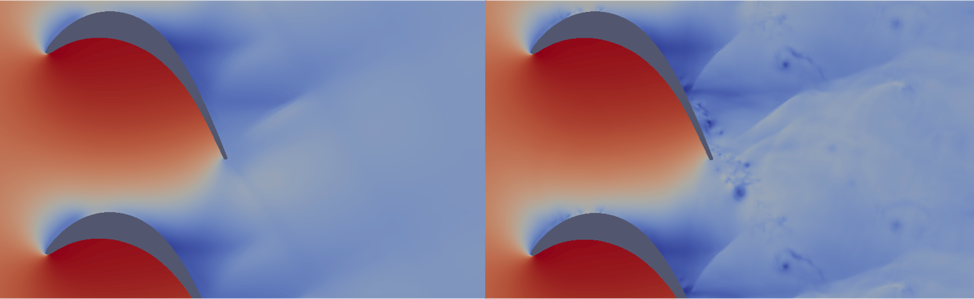}
\hfill \vspace{0.1mm} \includegraphics[width=1.0\textwidth]{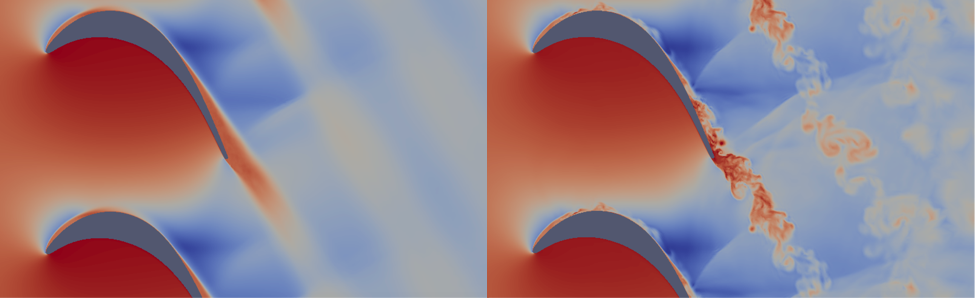}
\hfill \vspace{0.1mm} \includegraphics[width=1.0\textwidth]{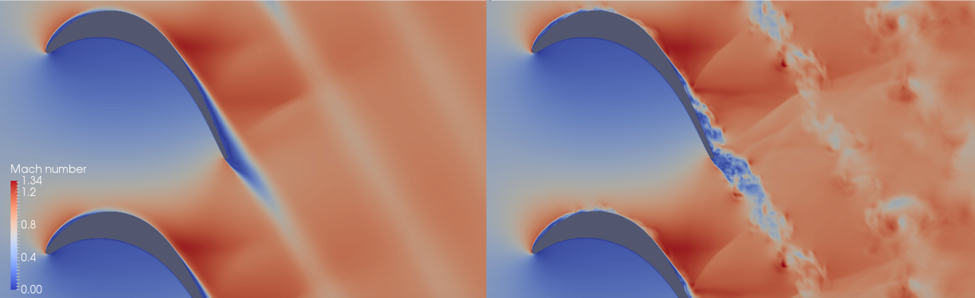}
\caption{\label{ptM} Pressure (top), temperature (center) and Mach number (bottom) fields on the periodic plane of the transonic T106C LPT. Time-averaged and instantaneous fields are shown on the left and right images, respectively. The unsteady shocks involved are resolved within one element.}
\end{figure}

We present ILES results for the three-dimensional transonic flow around the T106C low-pressure turbine (LPT) in off-design conditions \cite{Fernandez2018}. The isentropic Reynolds and Mach numbers on the outflow are $Re_{2,s} = 100,817$ and $M_{2,s} = 0.987$, respectively, and the angle between the inflow velocity and the longitudinal direction is $\alpha_{1} = 50.54 \ \textnormal{deg}$. Third-order HDG and DIRK(3,3) schemes are used for the discretization. The details of the simulation setup are presented in \cite{Fernandez2018}. Figure \ref{ptM} shows 2D slices of the time-averaged (left) and instantaneous (right) pressure, temperature and Mach number fields. Several unsteady shock waves that oscillate around a baseline position are present in this flow, as illustrated by the smoother shock profiles in the average fields compared to the instantaneous fields. These unsteady shocks are resolved within one element.

\subsection{Bibliography notes}

The HDG method for the incompressible Euler and Navier-Stokes equations was introduced in \cite{NguyenPeraireCockburn10AIAAINS,Nguyen2011h}, and further developed in \cite{Giorgiani2014,Lehrenfeld2016,NguyenPeraireCockburn10HDG,Nguyen2012,Rhebergen2012a,Ueckermann2016}. An analysis of the HDG method for the steady-state incompressible Navier-Stokes equations is presented in \cite{cesmelioglu2017analysis}. A superconvergent HDG method for the steady-state incompressible Navier-Stokes equations 
is developed in \cite{Qiu2016}. A comparison of HDG and finite volume methods for incompressible flows is presented in \cite{Ahnert2014}. No other schemes within the hybridized DG family, such as the EDG and the IEDG methods, have been applied to incompressible flows.

The HDG method for the compressible Euler and Navier-Stokes equations was first introduced in \cite{Peraire2010}, and further investigated in \cite{Nguyen2011b,Nguyen2012,Schutz2013,Schutz2013a,Woopen2014b,Fernandez2017a}. Additional developments of the HDG method for compressible flows include a multiscale method \cite{N2013}, a time-spectral method \cite{Chaurasia2013,Chaurasia2014}, and a viscous-inviscid monolithic solver \cite{Moro2013}. The Embedded Discontinuous Galerkin (EDG) and Interior Embedded Discontinuous Galerkin (IEDG) methods for the compressible Euler and Navier-Stokes equations were presented in \cite{Peraire2011} and \cite{Nguyen2015c,Fernandez2016}, respectively, and further investigated in \cite{Fernandez2017,Fernandez2017a}.

Other miscellaneous topics on hybridized DG methods for fluid flows include error estimation and adaptivity \cite{Balan2013,Dahm2014a,Fidkowski2016,Giorgiani2014,Jaust2014a,Moro2013,Moro2017,Woopen2014a,Woopen2014c}, entropy-stable formulations \cite{Fernandez:2018:entropyStable,williams2018entropy}, and shock capturing for steady \cite{Nguyen2011a,Moro2016} and unsteady \cite{Fernandez2018,Fernandez:shockCapturing:2018} flows. The relationship between the stabilization matrix and the resulting Riemann solver is investigated in \cite{Fernandez2017}. Finally, parallel implementation and efficiency considerations are discussed in \cite{Roca2013,Fernandez2017a}.
\section{Wave propagation in solids}

\subsection{Linear Elastodynamics}

Several HDG formulations have been proposed in the literature for linear elastic wave propagation. Each of them has pro and cons, and they will be briefly reviewed in subsection~\ref{subsec:bib_elastodyn}. We will only focus here on the \emph{velocity -- deformation-gradient} formulation, which is close to the HDG formulation we use for nonlinear elastodynamics.

\subsubsection{Governing equations}

We consider small transient adiabatic perturbations of an elastic body, which is at rest in a reference configuration $\Omega$. The perturbations are described using a deformation mapping $\bm{\varphi}$ between a reference configuration $\Omega$ and a current configuration $\Omega_t$ of the form $\bm{y} = \bm{\varphi}(\bm{X},t)$. Here, $\bm{X}$ is the coordinate in the reference configuration $\Omega$ and $\bm{y}$ denotes the position of material particle $\bm{X}$ after deformation at time $t$. The velocity is denoted by $\bm{v} = \partial_t \bm{\varphi}$, and the density of the reference configuration is denoted by $\rho$. Let $\bm{f}$ be the body force per unit reference volume. The motion of the elastic body under small perturbations is governed by the following linear elastic wave equation 
\begin{equation}
\label{LE}
\rho \, \partial_t \bm{v} + \nabla \cdot \bm{\sigma} = \bm{f}, \quad \mbox{in } \Omega \times (0,T) .
\end{equation}
where $\bm{\sigma}$ is the Cauchy stress tensor depending on two Lam\'{e} parameters $\lambda$ and $\mu$ for an isotropic body, and on the local state of deformation $\bm{\varphi}$. It is customary to write $\bm{\sigma}$ as a function of the infinitesimal strain tensor $\bm{\epsilon}$ under the assumption of small deformations. However, one could also directly write $\bm{\sigma}$ as a function of the deformation gradient $\bm{F}$, i.e.
\begin{equation}
\label{eq:sigma_cauchy}
    \bm{\sigma}(\bm{F}) = \mu \left( \bm{F} + \bm{F}^T \right) + \left( \lambda (\text{tr}(\bm{F}) - d ) -2\mu  \right) \bm{I}.
\end{equation}
Here $d$ is the spatial dimension of the problem, $\bm{I}$ is the identity tensor and $\bm{F}$ is the deformation gradient
\begin{equation}
\label{eq:deformation_grad}
    \bm{F} = \frac{\partial \bm{\varphi}}{\partial \bm{X}} .
\end{equation}

Now the elastic wave equation can be rewritten as
\begin{subequations} \label{LEM}
\begin{alignat}{3}
\partial_t \bm{F} - \nabla \bm{v} & = & 0, &  \quad \mbox{in } \Omega \times (0,T) , \label{eq:LEM1} \\
\rho \, \partial_t \bm{v} + \nabla \cdot \bm{\sigma}(\bm{F})  & = & \bm{f}, &  \quad \mbox{in } \Omega \times (0,T) , \label{eq:LEM2}
\end{alignat}
\end{subequations}
the first equation being the time derivative of~\eqref{eq:deformation_grad}.  The boundary conditions are given as
\begin{equation}
\label{NBCLE}
\begin{split}
\bm{v}  & = \bcD, \quad \mbox{on } \partial \Omega_D \times (0,T), \\
\bm{\sigma} \bm{n} & = \bcN, \quad \mbox{on } \partial \Omega_N \times (0,T),
\end{split}
\end{equation}
where $\partial \Omega_N$ is a part of the boundary $\partial \Omega$ such that $\partial \Omega_N \cup \partial \Omega_D = \partial \Omega$ and $\partial \Omega_N \cap \partial \Omega_D = \emptyset$.

\subsubsection{Formulation}

The HDG method seeks an approximation $(\bF,\bv,\bvh) \in \bm{\mathcal{Q}}_h^k \times \bm{\mathcal{V}}_h^k \times \bm{\mathcal{M}}_h^k$ such that

\begin{subequations}\label{eq:HDG_linelas}
\begin{alignat}{3}
  \Ohprod{\partial_t \bF}{\bm{G}} + \Ohprod{\bv}{\nabla \cdot \bm{G}} - \pOhprod{\bvh}{\bm{G} \bm{n}} =0 & & & \label{eq:HDG_linelas1}\\[1ex]
   \Ohprod{\rho \, \partial_t \bv}{ \bm{w}} + \Ohprod{\bm{\sigma}(\bF)}{\nabla \bm{w}}+\pOhprod{\bsh \bm{n}}{\bm{w}} -  \Ohprod{\bm f}{\bm{w}}= 0 & & &\label{eq:HDG_linelas2}\\[1ex]
 \pOhGprod{\bsh \bm{n}}{ \bm{\mu}} + \GprodN{ \bsh \bm{n} - \bcN}{\bm{\mu}} + \GprodD{ \bvh - \bcD}{\bm{\mu}}= 0 & & & \label{eq:HDG_linelas4} 
\end{alignat}
for all  $\bm{G} \in \bm{\mathcal{Q}}_{h}^k$, $\bm{w} \in \bm{\mathcal{V}}_h^k$, and $\bm{\mu}\in \bm{\mathcal{M}}_h^k$, where the numerical flux $\bsh$ is given by 
\begin{equation}
 \bsh \bm{n} = \bm{\sigma}(\bF) \bm{n} - \bm{S}( \bv - \bvh) . \label{eq:HDG_linelas5}
\end{equation}
\end{subequations}

Here we make use of the stabilization matrix $\bm{S}$ whose expression is discussed below. The first two equations \eqref{eq:HDG_linelas1}-\eqref{eq:HDG_linelas2} are obtained by multiplying the elastic wave equations~\eqref{eq:LEM1}-\eqref{eq:LEM2} by test functions and integrating the resulting equations by parts. The third equation~\eqref{eq:HDG_linelas4} enforces the continuity of the $L^2$ projection of the numerical flux $\bsh\bm{n}$ and imposes weakly the Dirichlet and Neumann boundary conditions. The last equation~\eqref{eq:HDG_linelas5} defines the numerical flux.

As a side note, the pointwise stress-strain relation~\eqref{eq:sigma_cauchy} could be applied through an element-based $L^2$-projection instead. For piecewise-constant Lam\'e parameters (used to produce the results shown in Table~\ref{tab:eoc_linear_elastodyn}) both approaches are equivalent and provide similar results.

\subsubsection{Stabilization matrix}


By using a simple dimensional analysis for the expression of the numerical flux \eqref{eq:HDG_linelas5}, it turns out that the stabilization matrix should be homogeneous to an impedance. It is a natural choice to consider impedances either the compressional elastic wave impedance, i.e. $\rho c_p$, or the shear wave impedance, i.e. $\rho c_s$. A simple choice for $\bm{S}$ would therefore be 
\begin{equation} \label{eq:stab_matrix_lin}
  \bm{S} =  \rho c_p \, \bm{I}, \ \ \ \ \text{or} \ \ \ \  \bm{S} =  \rho c_s \, \bm{I},
\end{equation}
where $c_p = \sqrt{(\lambda+2\mu) / \rho}$ is the compressional wave velocity, and $c_s = \sqrt{\mu / \rho}$ is the shear wave velocity. More sophisticated parameter-free $\bm{S}$ have been proposed in \cite{terrana2017spectral} do deal with impedance jumps at element boundaries, and acoustic waves coupling. However, for linear elastic wave problems, the accuracy of the approximation is only slightly dependent on $\bm{S}$, and therefore a wide range of values is acceptable for $\bm{S}$. It appears that, most of the time, choosing either one of the two impedances provides very satisfactory results.

\subsection{Nonlinear Elastodynamics}\label{subsec:nonlinear_elasdyn_HDG}

\subsubsection{Governing equations}

We now consider large time-dependent deformations of an elastic body defined by a deformation mapping $\bm{\varphi}$ between a reference configuration $\Omega$ and a current configuration $\Omega_t$ of the form $\bm{y} = \bm{\varphi}(\bm{X},t)$. Here, $\bm{X}$ is the coordinate in the reference configuration $\Omega$ and $\bm{y}$ denotes the position of material particle $\bm{X}$ after deformation at time $t$. The velocity is denoted by $\bm{v} = \partial_t \bm{\varphi}$, and the density of the reference configuration is denoted by $\rho$. Let $\bm{f}$ be the body force per unit reference volume. The boundary $\partial \Omega$ is divided into two complementary disjoint parts $\partial \Omega_D$ and $\partial \Omega_N$, where the prescribed deformation $\bcD$ and traction $\bcN$ are imposed, respectively. The motion of the elastic body under large deformations is governed by the following equations stated in Lagrangian form

\begin{subequations}
\label{UNEL}
\begin{alignat}{4}
\partial_t \bm{F}  - \nabla  \bm{v} & = & 0, & \quad \mbox{in } \Omega \times (0,T) , \label{eq:nlelastodyn_1} \\
 \rho \, \partial_t \bm{v} -\nabla \cdot \, {\bm{P}}  & = & \bm{f}, & \quad \mbox{in } \Omega \times (0,T),\label{eq:nlelastodyn_2} \\
\bm{P} - \bm{F} \SndPK (\bm{F}) & = & 0, & \quad \mbox{in } \Omega  \times (0,T), \label{eq:nlelastodyn_3}\\
\bm{v}  & = & \bcD, & \quad \mbox{on } \partial \Omega_D \ \times (0,T), \label{eq:nlelastodyn_4}\\
 {\bm{P}}  \bm{n} & = & \bcN, & \quad \mbox{on } \partial \Omega_N  \times (0,T),\label{eq:nlelastodyn_5}
\end{alignat}
\end{subequations}

The equation \eqref{eq:nlelastodyn_1} is just the time derivative of the definition of the gradient of deformation $\bm{F}$. The conservation of linear momentum and equation is stated with \eqref{eq:nlelastodyn_2}, and equation \eqref{eq:nlelastodyn_3} relates the first Piola-Kirchhoff tensor $\bm{P}$ with the second one $\SndPK$. The two last equations \eqref{eq:nlelastodyn_4}-\eqref{eq:nlelastodyn_5} express the boundary conditions. The gradient $\nabla$ and divergence $\nabla\cdot$ operators are taken with respect to the coordinate $\bm{X}$ of the reference configuration. To complete the problem description, an initial configuration $\bm{v}(\bm{X},t=0) = \bm{v}_0(\bm{X})$ and $\bm{F}(\bm{X},t=0) = \bm{F}_0(\bm{X})$ for all $\bm{X} \in \Omega$ has to be prescribed.

For hyperelastic materials the first and second Piola-Kirchhoff stress tensors $\bm{P}$ and $\SndPK$ are derived from a scalar strain energy function $\psi$ through
\begin{equation}\label{eq:rel_P_S_Psi}
  \bm P = \frac{\partial \psi (\bm{F})}{\partial \bm F} \quad \text{and} \quad \SndPK = 2 \frac{\partial \psi (\bm{F})}{\partial \bm C} ,
\end{equation}
with $\bm C = \bm F^T \bm F$ the right Cauchy-Green stress tensor. Hence, both $\bm P$ and $\bm S$ are functions of the deformation gradient and material parameters. For the applications in this section, a Saint Venant-Kirchhoff (SVK) model has been considered. For this model the second Piola-Kirchhoff tensors is given by
\begin{equation}\label{eq:psi_SVK_NeoHook}
  \SndPK =  \lambda \text{tr} (\bm{E})\bm{I} +2\mu\bm{E}
\end{equation}
whith the Lam\'e parameters $(\lambda,\mu)$, the Lagrangian strain tensor $\bm{E} = \frac{1}{2} (\bm{C} - \bm{I})$.

Below we introduce a HDG method for solving the nonlinear elasticity equations \eqref{UNEL}.

\subsubsection{Formulation}

We seek an approximation $(\bm{F}_h, \bm{P}_h, \bm{v}_h, \bvh) \in  \bm{\mathcal{Q}}_{h}^k \times \bm{\mathcal{Q}}_{h}^k \times \bm{\mathcal{V}}_h^k  \times \bm{\mathcal{M}}_h^k$ such that

\begin{subequations}\label{eq:HDG_nelas}
\begin{alignat}{3}
  \Ohprod{\partial_t \bF}{\bm{G}} + \Ohprod{\bv}{\nabla \cdot \bm{G}} - \pOhprod{\bvh}{\bm{G} \bm{n}} =0 & & & \label{eq:HDG_nelas1}\\[1ex]
   \Ohprod{\rho \, \partial_t \bv}{ \bm{w}} + \Ohprod{\bP}{\nabla \bm{w}}+\pOhprod{\widehat{\bP} \bm{n}}{\bm{w}} -  \Ohprod{\bm f}{\bm{w}}= 0 & & &\label{eq:HDG_nelas2}\\[1ex]
  \Ohprod{\bP}{\bm{Q}} - \Ohprod{\bF \SndPK(\bm F_h)}{\bm{Q}} = 0 & & &\label{eq:HDG_nelas3}\\[1ex]
 \pOhGprod{\widehat{\bP}\bm{n}}{ \bm{\mu}} + \GprodN{ \widehat{\bP} \bm{n} - \bcN}{\bm{\mu}} + \GprodD{ \bvh - \bcD}{\bm{\mu}}= 0 & & & \label{eq:HDG_nelas4} 
\end{alignat}
for all  $\bm{G}, \bm{Q} \in \bm{\mathcal{Q}}_{h}^k$ $\bm{w} \in \bm{\mathcal{V}}_h^k$, and $\bm{\mu}\in \bm{\mathcal{M}}_h^k$, where the numerical flux $\widehat{\bP}$ is given by 
\begin{equation}
 \widehat{\bP} \bm{n} = \bP \bm{n} - \bm{S}( \bv - \bvh) . \label{eq:HDG_nelas5}
\end{equation}
\end{subequations}

Here the stabilization tensor $\bm{S}$ does have an important effect on both the stability and accuracy of the method, and its design will be discussed below. Let us briefly comment on the equations defining the HDG method. The first two equations \eqref{eq:HDG_nelas1}-\eqref{eq:HDG_nelas2} are obtained by multiplying the governing equations~\eqref{eq:nlelastodyn_1}-\eqref{eq:nlelastodyn_2} by test functions and integrating the resulting equations by parts. The third equation~\eqref{eq:HDG_nelas3} is the weak version of~\eqref{eq:nlelastodyn_3}. The fourth equation~\eqref{eq:HDG_nelas4} enforces the continuity of the $L^2$ projection of the numerical flux $\widehat{\bP}\bm{n}$ and imposes weakly the Dirichlet and Neumann boundary conditions. The last equation~\eqref{eq:HDG_nelas5} defines the numerical flux.

\subsubsection{Stability}

We now give an insight into the energy evolution of our HDG method. Let us consider equation~\eqref{eq:HDG_nelas1} with test function $\bm{G} = \bP$, equation~\eqref{eq:HDG_nelas2} integrated by part with $\bm{w} = \bv$ and equation~\eqref{eq:HDG_nelas4} with $\bm{\mu} = \bvh$. After summing all these equations, and after some simplifications it comes the following energy identity
\begin{equation}\label{eq:energy_identity_nelas}
\Ohprod{\partial_t E_h}{1} +  \pOhprod{ \bm{S} (\bm{v}_{h} - \bvh) }{( \bm{v}_h - \bvh)} 
 =  \Ohprod{\bm{f}}{\bv} +  \GprodN{ \bm{t}}{\bvh} 
\end{equation}
with $\partial_t E_h$ the time derivative of the total discrete energy
\begin{equation}
  \partial_t E_h = \partial_t \left( \frac{1}{2}\rho \bv^2 \right) + \bP : \partial_t \bF
\end{equation}
where the last term is equal to $\partial_t \psi_h$, i.e. the time derivative of the discrete elastic potential energy. It comes from \eqref{eq:energy_identity_nelas} that, without external actions ($\bm{f} = 0$ and $\bcN = 0$), and if $\bm{S}$ is positive definite, the total energy decreases due to velocity jumps at element boundaries. The proposed HDG scheme is therefore stable with the jump term playing a stabilization role.

\subsubsection{Stabilization matrix}

The dimensional analysis done for the $\bm{S}$ in the linear case still applies for the nonlinear one. However, when large deformations occurs it becomes necessary to increase $\bm{S}$ in order to insure the convergence of the Newton's method. Therefore, one simple choice for the stabilization matrix is to scale the linear elastic $\bm{S}$ from~\eqref{eq:stab_matrix_lin} with an amplification factor, i.e.
\begin{equation} \label{eq:stab_matrix_nonlin}
  \bm{S} =  \alpha \, \rho c_p \, \bm{I}, \ \ \ \ \text{or} \ \ \ \  \bm{S} =  \alpha \, \rho c_s \, \bm{I},
\end{equation}
where the factor $\alpha$ is problem-dependent. In spite of its simplicity, the above stabilization tensor works well for many nonlinear test cases, with, most of the time, $\alpha \in [1,10]$. Contrary to the linear case, both the stability and the accuracy of our HDG scheme heavily depend on $\bm{S}$ and therefore the coefficient $\alpha$ plays a crucial role. Practically, $\alpha$ is determined after only a few trials.

We emphasize that it may not be a good idea to build $\bm{S}$ based on the material elasticity tensor, as proposed for elastostatics in~\cite{SoonPhdThesis,SoonCockburnStolarski09LE}. Indeed, in the linear elasticity case, this tensor is symmetric positive-definite everywhere in the domain. However, in the nonlinear elasticity case, it is generally no longer the case in regions where large deformation occurs. In this last configuration the energy identity~\eqref{eq:energy_identity_nelas} no longer holds and Newton's method typically fails to converge.



\subsection{Numerical Results}

\begin{figure}[!t]
  \centering
 \includegraphics[width=1.\linewidth]{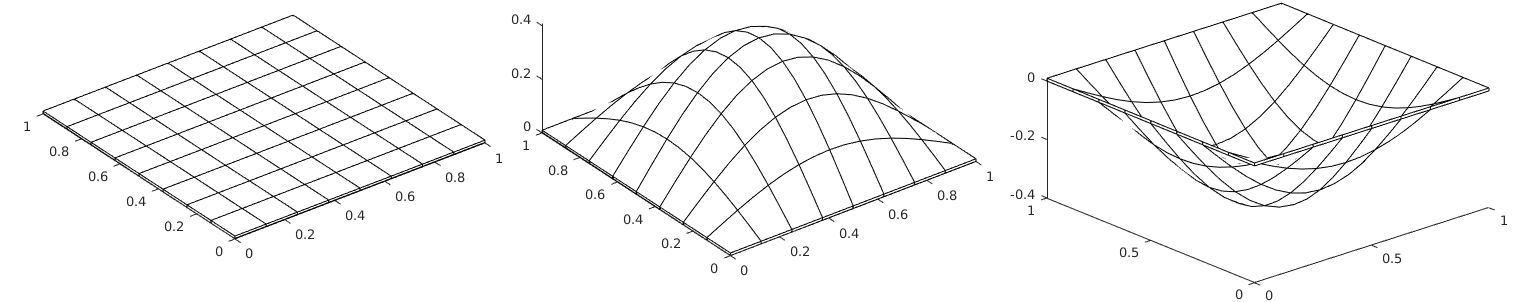}
    \caption{Clamped vibrating plate. Configurations at $t=0$ (left), $t=1$ (center) and $t=2$ (right).}
    \label{fig:vibrating_plate}
\end{figure}

We present here a simple numerical example in order to assess the convergence of our HDG formulations for both linear and nonlinear elastodynamics. In particular, we consider a square plate of dimensions $1 \times 1$ and of thickness $0.01$ that is clamped on its four sides, i.e.\ with homogeneous Dirichlet boundary conditions $\bcD=0$. The plate vibrates such that the exact deformation mapping -- illustrated in Figure \ref{fig:vibrating_plate} -- is
\begin{equation}
    \bm{\varphi}(\bm{X},t) = \left( 0,\ \  0, \ \  Z + 0.4 \sin(\pi t) \sin(\pi X) \sin(\pi Y) \right)^T
\end{equation}
where $\bm{X} = \left( X,Y,Z \right)^T$ are the positions in the undeformed initial configuration at $t=0$. Time dependent body forces, and tractions on the upper and lower surfaces are computed from the exact solution and imposed all along the simulation. The Lam\'e parameters are $\mu=1$ and $\lambda=1.5$, the density $\rho = 1$ and the stabilization matrix is set $\bm{S} = 2 \rho c_s \bm{I}$. The DIRK(3,3) scheme is used for the temporal discretization, and the time-step size is chosen sufficiently small so that the spatial discretization errors dominate. Both linear elastic and nonlinear hyperelastic (SVK) materials have been considered.

Numerical results are compared at $t=1$ with the exact ones for HDG-$\mathbb{P}_{k}$ with polynomial degrees $k\in \{1,2,3 \}$. The 3D mesh of the plate is uniformly refined in $\bm{e}_x$ and $\bm{e}_y$ directions. All simulations make use of only one element in the thickness direction. The $L^2$-errors with estimated orders of convergence (e.o.c) are reported in Table~\ref{tab:eoc_linear_elastodyn} for the linear case, and in Table~\ref{tab:eoc_nonlinear_elastodyn} for the nonlinear one.

It is worth noting that for the linear case, both the velocity and the gradient $\bF$ converge with the optimal order $k+1$. However, the analysis does not guarantee the optimal order of convergence for the gradient (see the discussion in subsection~\ref{subsec:bib_elastodyn}). For the nonlinear case, the velocity still converges optimally while the convergence order for the gradient $\bF$ is not clear, being almost $k+1/2$ for linear approximations, and somewhere beween $k+1/2$ and $k+1$ for $k\in \{ 2,3 \}$.  We emphasize that this last result depends on the choice of the stabilization matrix, and better convergence rates have been obtained when $\bm{S}$ is allowed to vary between simulations. However, an automatic optimal design of $\bm{S}$ for nonlinear elastic problems is still an open issue, and it is likely that $\bm{S}$ should be adaptive, as explained in \cite{EyckCelikerLew07part1,EyckCelikerLew08part2}.

\begin{table}
\centering
\begin{small}
\begin{tabular}{|c||c|c|c|c||c|c|c|c|}
  \hline
   &  \multicolumn{4}{c||}{ HDG-$\mathbb{P}_{1}$ } &  \multicolumn{4}{c|}{ HDG-$\mathbb{P}_{2}$ } \\
  $h$  & $\| \bm{v} - \bv \| $  & e.o.c & $\| \bm{F} - \bF \| $  & e.o.c & $\| \bm{v} - \bv \| $  & e.o.c & $\| \bm{F} - \bF \| $  & e.o.c \\
  \hline
  \hline
  0.5000 & 1.91E-2  & -     & 1.38E-1  & -     & 1.74E-03 &  -    & 1.46E-02 & -     \\  
  0.3333 & 9.96E-3  & 1.60  & 6.76E-2  & 1.76  & 4.60E-04 &  3.27 & 4.40E-03 & 2.97  \\
  0.2500 & 5.90E-3  & 1.82  & 3.92E-2  & 1.90  & 1.91E-04 &  3.06 & 1.56E-03 & 3.59  \\
  0.1666 & 2.71E-3  & 1.92  & 1.77E-2  & 1.95  & 6.19E-05 &  2.77 & 2.26E-04 & 4.77  \\
  0.1250 & 1.54E-3  & 1.96  & 1.01E-2  & 1.97  & 2.60E-05 &  3.01 & 9.30E-05 & 3.08  \\
  0.0909 & 8.22E-4  & 1.97  & 5.36E-3  & 1.96  & 9.71E-06 &  3.10 & 3.56E-05 & 3.02  \\
  0.0625 & 3.91E-4  & 1.98  & 2.55E-3  & 1.98  & 2.99E-06 &  3.14 & 1.00E-05 & 3.38  \\
  0.0435 & 1.90E-4  & 1.99  & 1.24E-3  & 1.99  & -        &   -   &      -   &   -   \\
  \hline
\end{tabular}

\vspace*{2.mm}

\begin{tabular}{|c||c|c|c|c|}
    \hline
    &  \multicolumn{4}{c|}{ HDG-$\mathbb{P}_{3}$} \\
   $h$  & $\| \bm{v} - \bv \| $  & e.o.c & $\| \bm{F} - \bF \| $  & e.o.c  \\
  \hline
  \hline
  0.5000 & 1.29E-04 & -     & 2.16E-03 &  -    \\
  0.3333 & 4.07E-05 & 2.84  & 4.14E-04 &  4.08 \\    
  0.2500 & 1.64E-05 & 3.17  & 7.41E-05 &  5.98 \\
  0.1666 & 2.86E-06 & 4.31  & 1.06E-05 &  4.80 \\
  0.1250 & 9.08E-07 & 3.98  & 3.23E-06 &  4.13 \\
  0.0909 & 2.52E-07 & 4.03  & 8.37E-07 &  4.24 \\
  \hline
  \end{tabular}
\end{small}
\caption{Numerical errors and estimated orders of convergence for the linear elastodynamic case.}
\label{tab:eoc_linear_elastodyn}
\end{table}


\begin{table}
\centering
\begin{small}
\begin{tabular}{|c||c|c|c|c||c|c|c|c|}
  \hline
   &  \multicolumn{4}{c||}{ HDG-$\mathbb{P}_{1}$ } &  \multicolumn{4}{c|}{ HDG-$\mathbb{P}_{2}$ } \\
  $h$  & $\| \bm{v} - \bv \| $  & e.o.c & $\| \bm{F} - \bF \| $  & e.o.c & $\| \bm{v} - \bv \| $  & e.o.c & $\| \bm{F} - \bF \| $  & e.o.c \\
  \hline
  \hline
  0.5000 & 6.25E-3  & -      & 4.06E-2  & -     & 2.37e-03 &  -     & 1.14e-02 &   -    \\  
  0.3333 & 7.00E-3  & -0.28  & 3.90E-2  & 0.10  & 8.47e-04 &  2.55  & 4.31e-03 &  2.40  \\
  0.2500 & 4.02E-3  & 1.93   & 2.75E-2  & 1.22  & 4.06e-04 &  2.56  & 1.99e-03 &  2.67  \\
  0.1666 & 1.65E-3  & 2.19   & 1.50E-2  & 1.49  & 1.47e-04 &  2.50  & 7.26e-03 &  2.49  \\
  0.1250 & 8.71E-4  & 2.22   & 9.91E-3  & 1.44  & 6.26e-05 &  2.97  & 3.41e-04 &  2.63  \\
  0.0909 & 3.96E-4  & 2.47   & 6.36E-3  & 1.40  & 2.60e-05 &  2.78  & 1.53e-04 &  2.52  \\
  0.0625 & 1.89E-4  & 1.97   & 3.81E-3  & 1.37  & 8.98e-06 &  2.84  & 5.84e-05 &  2.56  \\
  \hline
\end{tabular}

\vspace*{2.mm}

\begin{tabular}{|c||c|c|c|c|}
    \hline
    &  \multicolumn{4}{c|}{ HDG-$\mathbb{P}_{3}$} \\
   $h$  & $\| \bm{v} - \bv \| $  & e.o.c & $\| \bm{F} - \bF \| $  & e.o.c  \\
  \hline
  \hline
  0.5000 & 8.85e-04 &  -     & 3.49e-03 &    -   \\
  0.3333 & 2.78e-04 &  2.85  & 8.61e-04 &  3.45  \\    
  0.2500 & 1.18e-04 &  2.99  & 3.56e-04 &  3.06  \\
  0.1666 & 2.45e-05 &  3.86  & 6.86e-05 &  4.06  \\
  0.1250 & 8.65e-06 &  3.64  & 2.40e-05 &  3.66  \\
  0.0909 & 2.44e-06 &  3.97  & 7.04e-06 &  3.87  \\
  0.0625 & 5.46e-07 &  4.00  & 1.71e-06 &  3.78  \\
  \hline
  \end{tabular}
\end{small}
\caption{Numerical errors and estimated orders of convergence for the nonlinear (SVK) elastodynamic case.} 
\label{tab:eoc_nonlinear_elastodyn}
\end{table}


\subsection{Bibliography notes} \label{subsec:bib_elastodyn}

As explained in the introduction, two attractive features of the HDG methods are the optimal convergence of the approximate gradient and the superconvergence property. Namely, when a polynomial degree $k\geq 1$ is used to build the approximate primal solution and the approximate gradient, both of them may converge with an optimal order $k+1$, and the post-processed primal solution may then converge with an extra order $k+2$. Although the superconvergence property has been observed for numerous PDEs (see~\cite{Nguyen2012} among many others), it is not always guaranteed. This is especially true for linear elasticity since the symmetric nature of the infinitesimal strain and Cauchy stress tensors adds an extra difficulty. This difficulty has motivated the development of several HDG approaches, which are briefly reviewed here.

The first HDG method for linear elastostatics was introduced by~\cite{SoonPhdThesis,SoonCockburnStolarski09LE}, and makes use of a displacement-strain-stress formulation. Optimal convergence of the gradient and superconvergence were then numerically observed. However, the analysis \cite{fu2015analysis} of the same method demonstrated that although the displacement converges with order $k+1$, the symmetric part of the gradient converges with only $k+1/2$, and the antisymmetric part of the gradient with $k$. Moreover, numerical experiments illustrated that suboptimal convergence. Thus no superconvergence property is ensured with this method, although it may sometimes be observed. The strain-velocity formulation in~\cite{terrana2017spectral} extends the previous method to the elastic wave equations, with an emphasis on the proper design of $\bm{S}$ for heterogeneous media, following a methodology similar to~\cite{Bui2015}. Numerical results for the post-processed solution in~\cite{terrana2017spectral} confirmed the elastostatics results. The HDG formulation~\eqref{eq:HDG_linelas} presented in this paper, as well as the stress-velocity formulation presented in the frequency domain by~\cite{bonnasse2017hybridizable} are both variations of the initial~\cite{SoonCockburnStolarski09LE} HDG method (see~\cite{fu2015analysis}). Recently the $M$-decompositions theory has been used to modify that original HDG method such that it becomes superconvergent on 2D meshes by enriching the local gradient spaces (see~\cite{CockburnFu2018}).

Based on the study of superconvergent HDG methods for diffusion~\cite{cockburn2012conditions}, the method proposed in~\cite{Cockburn2013} makes use of approximate \emph{weakly symmetric stresses}. The superconvergence is then ensured, but at the cost of the extra computation of the approximate rotation tensor, and by enriching the gradient space with matrix bubble functions, which depend on the shape of the elements. However, it is difficult to extend this method to nonlinear elasticity since it involves the explicit inversion of the constitutive relation.

An alternative superconvergent HDG method was proposed for elastostatics~\cite{chen2016robust,qiu2018hdg}. In~\cite{Hungria2017}, a 3D time-harmonic elastodynamics version of this method was presented, with an analysis and some numerical experiments. These methods achieve an optimal $k+1$ convergence for the gradient, and ensure the superconvergence property at the cost of an extra polynomial degree $k+1$ for the the approximate displacement. However this computational overcost is small since the standard degree $k$ is still used for the approximate trace, so that the size of the global system is not increased.

Finally, a last family of HDG methods for linear elasticity has been presented in~\cite{Nguyen2011k,Nguyen2012}. Its elastodynamic version is based on a displacement gradient-velocity-pressure formulation, making use of the relation
\begin{equation}\label{eq:elas_stokes_like}
  \nabla \cdot \bm{\sigma} = \nabla \cdot \left( \mu \nabla \bm{u} - p \bm{I} \right) \ \ \ \ \text{with} \ \ \ \
  p = - (\lambda + \mu) \nabla \cdot \bm{u},
\end{equation}
in order to mimic the HDG formulation for the Stokes flow. Therefore it inherits all the superconvergence properties of HDG for the Stokes equations (see~\cite{Cockburn_HDG_Stokes,Cockburn2011,CockburnShiHDGStokes}). This method does not require any enrichment of the gradient space, and it allows for the treatment of nearly-incompressible elastic materials. However, this formulation has some drawbacks. The identity~\eqref{eq:elas_stokes_like} holds only for homogeneous $\mu$, and when normal stresses are applied as a boundary condition, the superconvergence may be lost (see~\cite{Nguyen2011h}).

Recently, a new Hybrid High Order (HHO) method was designed for linear elasticity in~\cite{di2015hybrid}. Contrary to HDG, the HHO methods are based on a primal formulation, i.e. the gradient is not considered as a separate variable. But like HDG, HHO methods make use of a static condensation procedure to solve a global system on the approximate traces, wich are typically polynomials of degree $k$, making the computational cost of both methods similar. Moreover, a locally reconstructed displacement field superconverges with a garanteed $k+2$ order of convergence. Interestingly, in~\cite{cockburn2016bridging}, the HHO method was recast into the HDG framework to study the hidden links between the two approaches.

The HDG literature for nonlinear elasticity is less abundant. A first HDG method was proposed in \cite{SoonPhdThesis} and later recast as a minimization  of a nonlinear functional in \cite{kabaria2015hybridizable}. Optimal convergence of the deformation and its gradient were numerically observed. The extension of this method to nonlinear elastodynamics is provided in this paper in subsection~\ref{subsec:nonlinear_elasdyn_HDG}.

In~\cite{Nguyen2012} a nonlinear elastodynamic HDG scheme was proposed using a deformation gradient-velocity-pressure formulation. Like for its linear counterpart, this formulation is attractive since it allows for the treatment of nearly-incompressible materials. The method presented in this paper is close, but it does not consider the pressure as a separate variable. Interestingly, both formulations seems to provide suboptimal convergence of the approximate gradient for $k=1,2$ while it is optimal for $k=3$.

Finally, an original Green strain-displacement-velocity formulation was proposed in~\cite{Sheldon2016} for the purpose of solving fluid-structure interaction problems. Observed convergence rates were $k+1$ for the approximate displacements and velocities, but only $k$ for the approximate strains.

\section{Electromagnetic wave propagation}

In this section, the HDG methods are extended to the generalized Maxwell's equations. 
The resolution of Maxwell's equations presents one major difference with the systems presented previously. Because of the presence of the vector operator $\bfcurl$ the electromagnetic field is  determined from its tangential component. As a consequence, the HDG methods are redefined with the introduction of tangential components. In addition, in low frequency regime, the Gauss's law needs to be numerically enforced on the model.  In case of the charge conservation is not satisfied on the discrete level, numerical errors and instabilities are introduced.  Many techniques have been developed to impose the charge conservation condition on the system matter \cite{Evans1988,Li2011}. Among them, the generalized Lagrange multiplier (GLM) method \cite{Dedner2002,Munz2000} enforces the divergence condition by solving a modified system where the constraint condition is imposed through the using of Lagrange multiplier. 
The nature of the correction allows the control of the propagation and the dissipation of divergence errors. This approach preserves the conservation form of the generalized Maxwell's equations at a minimal cost of introducing one additional scalar variable inside the system.

\subsection{Governing equations}

\subsubsection{Generalized Maxwell's equations}

The generalized Maxwell's equations are given by
\begin{subequations}
\label{eq:3Dcontin}
\begin{alignat}{4}
	\varepsilon\dt\bfE-\bfcurl\bfH & = & -\bfJ, &\text{ in } \Omega\times(0,T), \label{eq:max1} \\
    \mu\dt\bfH + \bfcurl\bfE & = & 0,& \text{ in } \Omega\times(0,T),\label{eq:max2}\\
    \nabla\cdot\bfE & = &\frac{\rho}{\varepsilon_0},& \text{ in } \Omega\times(0,T),\label{eq:max3} \\
    \nabla\cdot\bfH & = & 0,& \text{ in } \Omega\times(0,T),\label{eq:max4}
\end{alignat}
\end{subequations}
where $\bfE$ is the electric field, $\bfH$ the magnetic field, and $\bfJ$ the current density. In addition, $\varepsilon$, $\mu$ and $\rho$ denote the permittivity, permeability and the electric charge density, respectively. We assume boundary conditions of the form
\begin{equation}
\label{eq:boundary}
    \bfn\times\bfE\times\bfn  =  -\bfn\times\bfE^{\text{inc}}\times\bfn, 
    \text{ on } \partial \Omega \times (0,T) , 
\end{equation}
where $\bfn$ denotes the unit outward normal to $\partial\Omega$, and $(\bfE^{\text{inc}},\bfH^{\text{inc}})$ is the incident field. Finally, the system is supplemented with the initial conditions
\begin{subequations}
\begin{alignat}{4}
	\bfE & = &\bfE_0, &\text{ on } \Omega \times \{ t = 0 \}, \\
	\bfH & = &\bfH_0, &\text{ on } \Omega \times \{ t = 0 \},
\end{alignat}
\end{subequations}
where $\bfE_0$ and $\bfH_0$ are the initial electric and magnetic fields.

\subsubsection{Generalized Lagrange Multipliers}
In order to avoid instabilities and unphysical solutions related to electric field $\bfE$, we need to impose (\ref{eq:max3}) on the electromagnetic model.  The generalized Lagrange multiplier (GLM) method has been succesfully applied to Maxwell's equations \cite{Dedner2002,Munz2000}. The principle of the method is to introduce a new (non-physical) scalar field $\phi$ into the system (\ref{eq:max1})-(\ref{eq:max4}) through the differential operator $D(\phi)$ as follows
\begin{subequations}
\label{eq:GLM}
\begin{alignat}{4}
	\varepsilon\dt\bfE-\bfcurl\bfH + \nabla \phi & = & -
	\bfJ,& \text{ in } \Omega\times(0,T), \\
    \mu\dt\bfH+\bfcurl\bfE & = & 0,& \text{ in } \Omega\times(0,T),\\
    D(\phi) + \nabla\cdot\bfE & = & \frac{\rho}{\varepsilon_0},& \text{ in } \Omega\times(0,T).
\end{alignat}
\end{subequations}
In addition to (\ref{eq:boundary}), the following homogeneous Dirichlet condition is imposed
\begin{equation}
\phi = 0, \quad\text{on}\; \partial\Omega \times  (0,T) .
\end{equation}
In order to preserve hyperbolicity of the new system, the operator $D(\phi)$ is defined as follows
\begin{equation}
D(\phi) = \frac{1}{\alpha_1^2}\dt\phi + \frac{1}{\alpha_2^2}\phi , 
\end{equation}
where $\alpha_1\in\bbR^+$ and $\alpha_2\in\bbR^+$ are dimensionless coefficients that control the amount of artificial coupling between (\ref{eq:GLM}a)-(\ref{eq:GLM}b) and (\ref{eq:GLM}c). The resulting system is referred to as the generalized Lagrange multiplier formulation 
of the Maxwell's equations (GLM-Maxwell).


\subsection{Formulation}

To define the HDG method for the GLM-Maxwell equations we introduce the following space: 
\begin{equation}
\bm{{\mathcal{S}}}_{h}^k   = \big\{ \bm{\mu} \in [L^2(\mathcal{E}_h)]^m \  \ : \ (\bm{\mu} \circ \bm{\phi}^p_F) \in [\mathcal{P}^k(F_{ref})]^m \mbox{and } \bm \mu \cdot \bm n|_{\bm F} = 0, \, \ \forall F \in \mathcal{E}_h \big\} .
\end{equation}
Note that this space consists of vector-valued functions whose normal component, $\bm \mu^n := \bm n \, (\bm \mu \cdot \bm n)$, vanishes on every face $F$ of $\mathcal{E}_h$. In other words, we have $\bm \mu = \bm \mu^t := \bm n \times  \bm \mu \times\bm n$ for all $\bm \mu \in \bm{{\mathcal{S}}}_{h}^k$.

Multiplying the GLM-Maxwell equations by appropriate test functions and using the fact $\phi = 0$, the HDG discretization reads as follows: Find the approximate solution $\left(\bfE_h,\bfH_h,\phi_h,\hat{\bfE}^{t}_h \right) \in 
\bm{\mathcal{V}}_h^k \times \bm{\mathcal{V}}_h^k \times {{\mathcal{W}}}_h^k \times {\bm{\mathcal{S}}}_h^k$ such that the following equations are satisfied on each element $K$:
\begin{eqnarray}
  	\label{eq:DGLocal1}
    \left(\varepsilon\dt\bfE_h, \bfv\right)_{K} - 
    \left(\bfH_h, \bfcurl\bfv\right)_{K} - 
    \left<\hat{\bfH}^t_h,\bf v\times\bfn\right>_{\partial K} 
  -  \left(\ \phi_h, \nabla \cdot \bfv \right)_{K} & = &
    - \left(\bfJ,\bfv\right)_{K}, \label{eq:dmax1}\\
    \left(\mu\dt\bfH_h, \bfw\right)_{K} + 
    \left(\bfE_h, \bfcurl\bfw\right)_{K} + 
    \left<\hat{\bfE}^t_h,\bfw\times\bfn\right>_{\partial K} & = &
    0, \label{eq:dmax2}\\
    \frac{1}{\alpha_1^2} \left(\dt \phi_h,\psi\right)_{K} + 
    \frac{1}{\alpha_2^2}\left(\phi_h,\psi\right)_{K} +   \frac{1}{\alpha_3^2}\left<\phi_h,\psi\right>_{\partial K} +
    \left( \nabla\cdot \bfE_h, \psi \right)_{K} & = &
    \left(\frac{\rho}{\varepsilon_0},\psi\right)_{K} .\label{eq:dmax3}
\end{eqnarray}
Note that $\hat{\bfH}^{t}_h$ and $\hat{\bfE}^{t}_h$ denote the approximate trace of ${\bfH}^{t} := \bm n \times \bm h \times \bm n$ and ${\bfE}^{t} := \bm n \times \bm e \times \bm n$ on the element boundaries, respectively. Next, we define $\hat{\bfH}^{t}_h$ as follows
\begin{equation}
\hat{\bfH}^{t}_h  =  \bfH^{t}_h -  \tau\left(\bfE_h  - \hat{\bfE}^{t}_h \right)\times \bfn \label{eq:hatH}
\end{equation}
where $\tau$ is a local stabilization parameter, and enforce the conservativity condition and the boundary condition as follows
\begin{equation}
\left< \bm n \times \hat{\bfH}^{t}_h, \bm \mu \right>_{\partial \mathcal{T}_h \backslash \partial \Omega} + \left< \hat{\bfE}^{t}_h + \bfn\times\bfE^{\text{inc}}\times\bfn, \bm \mu \right>_{\partial \Omega} = 0 .
\label{eq:hatBC}
\end{equation}
The test functions are taken as $\left(\bm v,\bm w,\psi,\bm \mu \right) \in 
\bm{\mathcal{V}}_h^k \times \bm{\mathcal{V}}_h^k \times {{\mathcal{W}}}_h^k \times {\bm{\mathcal{S}}}_h^k$. The term $\frac{1}{\alpha_3^2}\left<\phi_h,\psi\right>_{K} $ is added in \eqref{eq:dmax3} to provide additional stabilization of the divergence-free constraint.

According to the discussion on time-marching techniques in Section \ref{sec:timemethods}, the semi-discrete HDG formulation can be written as the DAE system (\ref{ch5appendx1:eq8}) with $\bm u$ being the vector degrees of freedom of $\left(\bfE_h,\bfH_h, \phi_h\right)$ and $\bm v$ being the vector of degrees of freedom of $\hat{\bm e}_h^t$. Note also that when constructing the global linear system for the degrees of freedom of $\hat{\bm e}_h^t$, we locally eliminate the degrees of freedom of $\phi_h$ by substituting it from (\ref{eq:dmax3}) into (\ref{eq:DGLocal1}). Therefore,  the introduction of the Lagrange multiplier $\phi_h$ does not affect the computational complexity of the HDG method. In other words, the computational complexity of the proposed HDG method is the same as that of the HDG methods presented in \cite{Christophe2018,Li2014,Li2015,Nguyen2011j}. Unlike the proposed HDG method, those HDG methods do not discretize the divergence-free constraint.

\subsection{Stability and consistency}


We can show that the local problem \eqref{eq:DGLocal1}-\eqref{eq:dmax3} is well-defined. Indeed, inserting (\ref{eq:hatH}) into (\ref{eq:DGLocal1}) and summing up the three equations (\ref{eq:DGLocal1})-(\ref{eq:dmax3}) yields
\begin{multline}
    \left(\varepsilon\dt\bfE_h, \bfv\right)_{K}  +     \left(\mu\dt\bfH_h, \bfw\right)_{K}    + \frac{1}{\alpha_1^2} \left(\dt \phi_h,\psi\right)_{K} + 
    \frac{1}{\alpha_2^2}\left(\phi_h,\psi\right)_{K} +   \frac{1}{\alpha_3^2}\left<\phi_h,\psi\right>_{\partial K} +     \left<\tau \bfE_h\times \bm n ,\bf v\times\bfn\right>_{\partial K} \\
    =     \left<\tau \hat{\bfE}_h^t\times \bm n ,\bf v\times\bfn\right>_{\partial K}   - \left(\bfJ,\bfv\right)_{K}   -     \left<\hat{\bfE}^t_h,\bfw\times\bfn\right>_{\partial K}  +  \left(\frac{\rho}{\varepsilon_0},\psi\right)_{K} .
\end{multline}
Integrating this equation from time $t_1=0$ to $t_2 = t$ and choosing $\left(\bm v,\bm w,\psi \right) = \left(\bfE_h,\bfH_h,\phi_h\right)$ as test functions, we obtain
\begin{multline}
\frac{1}{2}    \left(\varepsilon\bfE_h(t), \bfE_h(t) \right)_{K}  +  \frac{1}{2}   \left(\mu \bfH_h(t), \bfH_h(t) \right)_{K}    + \frac{1}{2\alpha_1^2} \left(\phi_h(t),\phi_h(t)\right)_{K} + \\
\int_{0}^t   \left( \frac{1}{\alpha_2^2}\left(\phi_h,\phi_h \right)_{K} +   \frac{1}{\alpha_3^2}\left<\phi_h,\phi_h\right>_{\partial K} +     \left<\tau \bfE_h\times \bm n ,\bm e_h \times\bfn\right>_{\partial K} \right) ds \\
= \frac{1}{2}    \left(\varepsilon\bfE_h(0), \bfE_h(0) \right)_{K}  +  \frac{1}{2}   \left(\mu \bfH_h(0), \bfH_h(0) \right)_{K}    + \frac{1}{2\alpha_1^2} \left(\phi_h(0),\phi_h(0)\right)_{K} + \\
        \int_{0}^t   \left( \left<\tau \hat{\bfE}_h^t\times \bm n ,\bm e_h \times\bfn\right>_{\partial K}   - \left(\bfJ,\bm e_h\right)_{K}   -     \left<\hat{\bfE}^t_h,\bm h_h \times\bfn\right>_{\partial K}  +  \left(\frac{\rho}{\varepsilon_0},\phi_h\right)_{K} \right) ds .
\end{multline}
This identity implies the local problem has a unique solution. 

In a similar manner, the following energy identity holds for the semi-discrete HDG formulation:
\begin{multline}
\frac{1}{2}    \left(\varepsilon\bfE_h(t), \bfE_h(t) \right)_{\calT_{h}}  +  \frac{1}{2}   \left(\mu \bfH_h(t), \bfH_h(t) \right)_{\calT_{h}}    + \frac{1}{2\alpha_1^2} \left(\phi_h(t),\phi_h(t)\right)_{\calT_{h}} + \\
\int_{0}^t   \left( \frac{1}{\alpha_2^2}\left(\phi_h,\phi_h \right)_{\calT_{h}} +   \frac{1}{\alpha_3^2}\left<\phi_h,\phi_h\right>_{\partial \calT_{h}} +     \left<\tau (\bfE_h - \hat{\bfE}_h^t)\times \bm n ,(\bm e_h - \hat{\bfE}_h^t) \times\bfn\right>_{\partial \calT_{h}}  \right) ds \\
    =   \frac{1}{2}    \left(\varepsilon\bfE_h(0), \bfE_h(0) \right)_{\calT_{h}}  +  \frac{1}{2}   \left(\mu \bfH_h(0), \bfH_h(0) \right)_{\calT_{h}}    + \frac{1}{2\alpha_1^2} \left(\phi_h(0),\phi_h(0)\right)_{\calT_{h}} + \\
 \int_{0}^t   \left(  - \left(\bfJ,\bm e_h\right)_{K}   +  \left(\frac{\rho}{\varepsilon_0},\phi_h\right)_{K} - \left<  \bm n \times \hat{\bfH}^{t}_h,  \bfn\times\bfE^{\text{inc}}\times\bfn \right>_{\partial \Omega}  \right) ds .
\end{multline}
This energy identity shows the existence and uniqueness of the numerical solution. In addition, the discrete energy
\begin{equation}
E(t) = \frac{1}{2}    \left(\varepsilon\bfE_h(t), \bfE_h(t) \right)_{\calT_{h}}  +  \frac{1}{2}   \left(\mu \bfH_h(t), \bfH_h(t) \right)_{\calT_{h}}    + \frac{1}{2\alpha_1^2} \left(\phi_h(t),\phi_h(t)\right)_{\calT_{h}} 
\end{equation}
decays in time whenever $\bm j = 0$, $\rho = 0$, and $\bfE^{\text{inc}} = 0$. Hence, the HDG method is well-defined and stable. Finally, it is easy to show that the exact solution also satisfies the HDG formulation (\ref{eq:DGLocal1})-(\ref{eq:hatBC}). Therefore, the HDG method is consistent.

\subsection{Numerical results}
In order to demonstrate the convergence and accuracy of the HDG method, a three-dimensional problem with no electric charge density (i.e.\ $\rho =0$) is considered. This problem involves the propagation of a standing wave in a cubic cavity $\Omega = (0,1)\times(0,1)\times(0,1)$ with perfect electrical conductor (PEC) boundaries up to a final time $T=1$. The permittivity is $\varepsilon_r = 2$, the permeability $\mu_r = 1$, and the current density is neglected, i.e.\ $\bfJ = 0$. The exact solution of the problem is given by
    \begin{equation*}
      \label{eq:cavityexact}
      \bfE(\bfx,t) = 
		\begin{bmatrix}
    		  \sin(\omega y)\sin(\omega z)\sin(\omega t)\\
    		  \sin(\omega x)\sin(\omega z)\sin(\omega t)\\
    		  \sin(\omega y)\sin(\omega x)\sin(\omega t)
		\end{bmatrix}, \quad
		\bfH(\bfx,t) = 
		\begin{bmatrix}
    		\left(\cos(\omega y)\sin(\omega x) - \cos(\omega z)\sin(\omega x)\right)\cos(\omega t)\\
    		\left(\cos(\omega z)\sin(\omega y) - \cos(\omega x)\sin(\omega y)\right)\cos(\omega t)\\
    		\left(\cos(\omega x)\sin(\omega z) - \cos(\omega y)\sin(\omega z)\right)\cos(\omega t)\\
		\end{bmatrix} , 
    \end{equation*}
where the angular frequency (or pulsation) is $\omega = 1$. The GLM coefficients are set to $\alpha_1 = \alpha_2 = \alpha_3 = 1$, and the stabilization parameter to $\tau = 2$. The DIRK(3,3) scheme is used for the temporal discretization, and the time-step size is chosen sufficiently small so that the spatial discretization errors dominate.

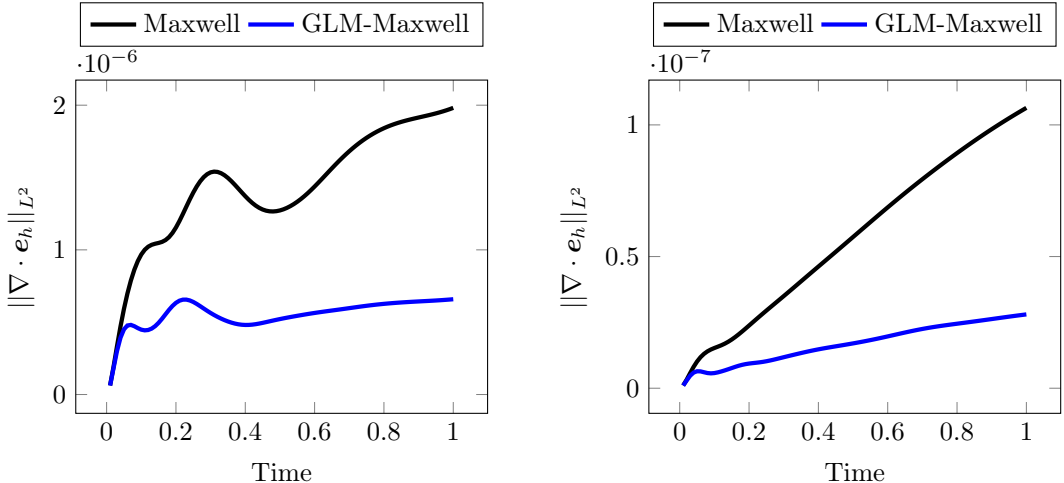
\begin{figure}[!h]
    \begin{tikzpicture}
     \begin{axis}[
      height=6cm,
      width=7cm, 
      ylabel style={at={(axis description cs:.15,.5)},anchor=south},
      xlabel={Time}, 
      ylabel={$||\nabla\cdot \bfE_h||_{L^2}$},
      legend style={at={(0.5,1.1)},anchor=south},
      legend entries={Maxwell,GLM-Maxwell},
      legend columns=2]    
    \addplot[mark=none,draw=black,ultra thick] table[x=time,y=div]{divP2n4.txt};
    \addplot[mark=none,draw=blue,ultra thick] table[x=time,y=div]{divfreeP2n4.txt};
     \end{axis}
     \end{tikzpicture}
     \hspace*{6mm}
    \begin{tikzpicture}
     \begin{axis}[
      height=6cm,
      width=7cm, 
      ylabel style={at={(axis description cs:.1,.5)},anchor=south},
      xlabel={Time}, 
      ylabel={$||\nabla\cdot \bfE_h||_{L^2}$},
      legend style={at={(0.5,1.1)},anchor=south},
      legend entries={Maxwell,GLM-Maxwell},
      legend columns=2]   
    \addplot[mark=none,draw=black,ultra thick] table[x=time,y=div]{divP3n4.txt};
    \addplot[mark=none,draw=blue,ultra thick] table[x=time,y=div]{divfreeP3n4.txt};
     \end{axis}
     \end{tikzpicture}
    \caption{Comparison of the numerical errors in the $L^2$ norm of $\nabla\cdot \bfE_h$ computed with the uncorrected Maxwell's equations (black) and the corrected GLM-Maxwell model (blue): HDG-$\mathbb{P}_{2}$ (left) and HDG-$\mathbb{P}_{3}$ (right).}
    \label{fig:divEP2P3}
 \end{figure}

Tables \ref{tab:tab1}, \ref{tab:tab2} and \ref{tab:tab3} present the numerical errors and estimated orders of convergence (e.o.c) for HDG-$\mathbb{P}_{k}$ {with polynomial degrees $k = 1$, $2$ and $3$}, respectively.  The convergence rates and errors for $\phi_h$ are shown in Table \ref{tab:tab4}. We observe that the convergence rates are optimal for all the variables. Finally, we compare the time evolution of the $L^2$-error norm of $\nabla\cdot \bfE_h$ for the uncorrected Maxwell's equations and the GLM-Maxwell system. In particular, Figure \ref{fig:divEP2P3} shows the time evolution for various polynomial orders on a $8 \times 8 \times 8$ mesh. We observe that the errors with the GLM-Maxwell model are smaller than those with the uncorrected Maxwell's equations. Therefore, the numerical treatment of the divergence-free constraint using the GLM-Maxwell model enhances accuracy and long-time stability.


\begin{table}[!h] 
\centering
\begin{scriptsize}
\begin{tabular}{| c || c |  c |  c | c  || c | c | c | c  |}
\hline
  &  \multicolumn{4}{c||}{ $L^2$ norm } &  \multicolumn{4}{c|}{ $H(\bfcurl)$ norm } \\
$h$  &  $\| \bfH - \bfH_h \| $  & e.o.c & $\| \bfE - \bfE_h \| $  & e.o.c & $\| \bfH - \bfH_h \| $  & e.o.c & $\| \bfE - \bfE_h \| $  & e.o.c \\
\hline	
\hline	
  $1/2 $ & 4.25E-02  & {-} & 7.41E-02  &{-} &  
  1.98E-01	& {-} & 6.59E-01 & {-} \\
  $1/4$ & 8.94E-03	& {2.2} & 1.08E-02  &{2.8} & 
   4.72E-02 & {2.0} & 2.34E-01 & {1.5} \\
  $1/8$ & 1.88E-03	& {2.3} & 1.97E-03  &{2.5} & 
  1.97E-02	& {1.3} & 9.67E-02 & {1.3} \\
\hline
\end{tabular}
\caption{Numerical errors and estimated orders of convergence with HDG-$\mathbb{P}_{1}$.}
\label{tab:tab1}
\end{scriptsize}
\vspace*{4mm}
\begin{scriptsize}
\begin{tabular}{| c || c |  c |  c | c  || c | c | c | c  |}
\hline	
  &  \multicolumn{4}{c||}{ $L^2$ norm } &  \multicolumn{4}{c|}{ $H(\bfcurl)$ norm } \\
$h$  &  $\| \bfH - \bfH_h \| $  & e.o.c & $\| \bfE - \bfE_h \| $  & e.o.c & $\| \bfH - \bfH_h \| $  & e.o.c & $\| \bfE - \bfE_h \| $  & e.o.c \\
\hline	
\hline	
  $1/2 $ & 9.15E-03     & {-}   & 1.86E-02  & {-}   & 
  5.40E-02 & {-}   & 1.45E-01 & {-}  \\
  $1/4$ & 3.23E-04	& {4.8} & 6.76E-04  &{4.8} & 
  7.56E-03 & {2.8} & 1.58E-02 & {3.2} \\
  $1/8$ & 2.75E-05	& {3.6} & 5.13E-05  &{3.7} & 
  1.54E-03	& {2.3} & 2.85E-03 & {2.5} \\
\hline
\end{tabular}
\caption{Numerical errors and estimated orders of convergence with HDG-$\mathbb{P}_{2}$.}
\label{tab:tab2}
\end{scriptsize}
\vspace*{4mm}
\begin{scriptsize}
\begin{tabular}{| c || c |  c |  c | c  || c | c | c | c  |}
\hline	
  &  \multicolumn{4}{c||}{ $L^2$ norm } &  \multicolumn{4}{c|}{ $H(\bfcurl)$ norm } \\
$h$  &  $\| \bfH - \bfH_h \| $  & e.o.c & $\| \bfE - \bfE_h \| $  & e.o.c & $\| \bfH - \bfH_h \| $  & e.o.c & $\| \bfE - \bfE_h \| $  & e.o.c \\
\hline	
\hline	
  $1/2 $ & 1.06E-03     & -   & 1.04E-03  & - & 
  9.14E-03 & -   & 1.39E-02 & -   \\
  $1/4$ & 1.50E-05	& 6.1 & 1.33E-05  & 6.3 & 
  3.36E-04 & 4.8 & 4.73E-04 & 4.9 \\
  $1/8$ & 5.73E-07	& 4.7 & 4.54E-07  & 4.9 & 
  2.69E-05	& 3.7 & 3.55E-05 & 3.7 \\
\hline
\end{tabular}
\caption{Numerical errors and estimated orders of convergence with HDG-$\mathbb{P}_{3}$.}
\label{tab:tab3}
\end{scriptsize}
\end{table} 

\begin{table}[!h] 
\centering
\begin{scriptsize}
\begin{tabular}{| c || c |  c |  c | c | c |  c |}
\hline	
  & \multicolumn{2}{c|}{HDG-$\mathbb{P}_{1}$} & \multicolumn{2}{c|}{HDG-$\mathbb{P}_{2}$}  & \multicolumn{2}{c|}{HDG-$\mathbb{P}_{3}$} \\
$h$  &  $\| \phi - \phi_h \| $  & e.o.c & $\| \phi - \phi_h \| $  & e.o.c & $\| \phi - \phi_h \| $  & e.o.c \\
\hline	
\hline	
  $1/2 $ & 8.90E-03 & - & 7.39E-04  & - & 1.48E-04  & - \\
  $1/4$ & 1.13E-03	& 3.0 & 2.52E-05  & 4.9 & 1.80E-06  & {6.4} \\
  $1/8$ & 2.07E-04	& 2.5 & 1.59E-06  & 4.0 & 5.86E-08  & {4.9} \\
\hline
\end{tabular}
\caption{Numerical errors and estimated orders of convergence for the Lagrange multiplier $\phi$ with HDG-$\mathbb{P}_{k} , \ k = \{ 1 , 2 , 3 \}$.}
\label{tab:tab4}
\end{scriptsize}
\end{table}

 \subsection{Bibliography notes}

Discretization of Maxwell's equations in high-frequency regime using HDG methods has been done in both frequency and time domains  \cite{Chen2017,Christophe2018,feng2016,li:2017,Li2014,Li2015,lu2017absolutely,Nguyen2011j,Vidal-Codina2018,zhu2017}. The first HDG method for the time-harmonic Maxwell's equations was proposed in \cite{Nguyen2011j} for two-dimensional problems. The extension
of the method to three dimensions was presented in \cite{Li2014,Li2015}. HDG was employed for full 3D modeling of the resonant transmission of THz waves through annular gaps in the field of nanoplasmonics \cite{Park2014,Yoo2016}.  An HDG method for computing nonlocal electromagnetic effects in three-dimensional metallic nanostructures has been recently introduced in \cite{Vidal-Codina2018}.  


\section{Perspectives}


In spite of considerable effort towards making DG methods more robust and computationally efficient, there are still open problems demanding advances on several research fronts. We end this paper with perspectives on ongoing extension and new development of hybridized DG methods for wave propagation problems. While HDG has been only applied to a wide variety of wave propagation problems, EDG and IEDG have been applied to compressible flows. Due to their significantly lower computational cost, the application of EDG and IEDG methods to solid mechanics, incompressible flows, and electromagnetism is encouraged.


Hybridized DG methods use polynomial spaces to approximate the solution on elements and faces. A possible extension is the enrichment of the approximation spaces with non-polynomial functions in order to capture discontinuities, singularities, and boundary layers. The hybridized DG framework may lend itself for this task because the enrichment can be done at the element level thanks to the discontinuous nature of the approximation spaces. Indeed, an HDG method using exponential kernels for high-frequency wave propagation is proposed in \cite{NguyenN.C.2015}, and an extended HDG method with heaviside enrichment for heat bimaterial problems is developed in \cite{Gurkan2017}.

 
In this paper, we have exclusively focused on implicit hybridized DG methods. It is highly desirable to develop hybridized DG methods that can be coupled with explicit time discretization for time-dependent problems. They should be computationally competitive to other explicit DG methods, while retaining some important advantages such as the superconvergence property.  As a step in this direction, explicit HDG methods have been devised for the acoustics wave equation \cite{StanglmeierNguyenPeraireCockburn16}. While extension of the explicit HDG methods to elastodynamics and electromagnetics is quite straightforward, it is not trivial to develop efficient explicit HDG methods for fluid dynamics. Another area of interest is to devise hybridized DG methods coupled with implicit-explicit (IMEX) time-marching schemes. This is recently persued for acoustics wave problems \cite{Kolkman2018}.

Also, the time-marching schemes for hybridized DG methods discussed herein are dissipative in the sense that the discrete energy is decaying in time for problems in which the exact energy is invariant in time. For many wave propagation problems, it is crucial to equip numerical methods with desirable conservation properties such as energy and momentum conservations for long-time simulations. There have been recent work on symplectic HDG methods for acoustic waves  \cite{Cockburn2017,SANCHEZ2017951}. It will be interesting to develop symplectic HDG methods for shallow water waves, elastic waves, and electromagnetic waves. 

Finally, we point out other work on the development of HDG methods for wave propagation problems. The first HDG method for the Helmholtz equation was introduced in \cite{GriesmaierMonk}. In \cite{Feng2012}, a wide family of discontinuous Galerkin methods, which included the HDG methods, were proven to be stable regardless of the wave number. The methods used piecewise linear approximations. In \cite{Cui2014}, an analysis of the HDG methods for the Helmholtz equations was carried which shows that the method is stable for any wave number, mesh and polynomial degree and which recovers the orders of convergence and superconvergence obtained previously in \cite{GriesmaierMonk}. In \cite{gopalakrishnan2015spectral}, the HDG method for eigenvalue problems was developed and analyzed.  HDG methods for the Oseen equations were developed and analyzed in  \cite{CesmeliogluCockburnHDGOseen13}. A systematic way of defining HDG methods for Friedrichs' systems has been developed in \cite{Bui2015}.  An explicit HDG method for Serre-Green-Naghdi wave model is devised in \cite{Samii2018}. The first HDG method for solving Korteweg-de Vries (KdV) type equations is developed and analyzed in \cite{Dong2017}.  An HDG method for coupled fluid-structure interaction problems is presented in \cite{Sheldon2016}. Hybridized DG methods for ideal and resistive MHD problems are recently developed in \cite{Ciuca2018}.



\section*{Acknowledgements}
The authors acknowledge the Air Force Office of Scientific Research (FA9550-15-1-0276 and FA9550-16-1-0214), the NASA (NNX16AP15A), and Pratt \& Whitney for supporting this work. P. Fernandez also acknowledges the financial support from the Zakhartchenko and ``la Caixa'' Fellowships.

\bibliographystyle{plain}
\bibliography{../bibliography/library,../bibliography/library_sebastien,../bibliography/library_pablo}
\end{document}